\documentclass[10pt]{amsart}
\usepackage{amsthm,amsmath,amssymb,amscd}
\usepackage{amssymb}
\usepackage{graphics}

\newtheorem{teor}{Theorem}[section]

\newtheorem{prop}{Proposition}[section]
\newtheorem{corol}{Corollary}[section]
\newtheorem{lemma}{Lemma}[section]

\theoremstyle{definition}
\newtheorem{dfn}[subsection]{Definition}

\newtheorem{rmk}[subsection]{Remark}

\newtheorem{defin}[subsection]{Definition}

\newcommand{\K}{K\"{a}hler}

\newcommand{\R}{\mathbb{R}}

\newcommand{\C}{\mathbb{C}}

\newcommand{\e}{\varepsilon}
\newcommand{\del}{\partial}

\newcommand{\frh}{{\mathfrak h}}
\newcommand{\frk}{{\mathfrak k}}

\newcommand{\ci}{{\mathcal C}^{\infty}}
\newcommand{\ip}[1]{\langle#1\rangle}
\newcommand{\rd}{\mathrm d}

\begin{document}

\title{Extremal metrics on blow ups}

\author{C. Arezzo} \author{F. Pacard } \author{M. Singer}

\maketitle
\section{Introduction}

In this paper we study the problem of constructing extremal \K\
metrics on blow ups at finitely many points of \K\ manifolds which
already carry an extremal metric.

In \cite{ca}, \cite{ca2} Calabi has proposed, as best
representatives of a given \K\ class $[\omega]$ of a complex compact
manifold $(M,J)$, a special type of metrics baptized {\em extremal}.
These metrics are critical points of the $L^{2}$-square norm of the
scalar curvature ${\bf s}$. The corresponding Euler-Lagrange
equation reduces to the fact that
\[
\Xi_{\bf s} : = J \, \nabla {\bf s} + i \, \nabla {\bf s}
\]
is a holomorphic vector field on $M$. In particular, the set of
extremal metrics contains the set of constant scalar curvature \K \
ones. Calabi's intuition of looking at extremal metrics as canonical
representatives of a given \K\ class has found a number of important
confirmations and also (unfortunately) nontrivial constraints.
Calabi himself proved that an extremal K\"ahler metric must have the
maximal possible symmetry allowed by the complex manifold $M$, and,
as observed by LeBrun and  Simanca \cite{ls}, this symmetry group
can be fixed in advance. More precisely, the identity component of
the isometry group of any extremal metric $g$ must be a maximal
compact subgroup of $\mbox{Aut}_0(M,J)$, the identity component of
the group $\mbox{Aut} (M,J)$ of biholomorphic maps of $M$ to itself.
This group thus contains the complexification of the isometry group,
but may be strictly larger (the blow-up of ${\mathbb P}^{2}$ at a
point is the simplest example of such a situation). Moreover Lebrun
and Simanca \cite{ls2} have proved that the set of \K\ classes
having an extremal representative is an open subset of $H^{1,1}(M,
{\mathbb C}) \cap H^2 (M, {\mathbb R})$ and Chen and Tian \cite{ct}
have proved the uniqueness of such metrics in a given \K\ class up
to automorphisms. Also, the important relationship between the
existence of extremal metrics and various stability notions of the
corresponding polarized manifolds (algebraic if the class is
rational, analytic otherwise) has been deeply investigated for
example by Tian \cite{ti}, Mabuchi \cite{ma} and Szekelyhidi
\cite{sz}. Yet, a complete understanding of the existence theory for
extremal metrics is still missing. Given this last fact, the first
two authors have started in \cite{ap} and \cite{ap2} to develop a
perturbation theory for constant scalar curvature \K\ metrics,
giving sufficient conditions for the existence of constant scalar
curvature \K\ metrics on the blow up at finitely many points of a
manifold which already carries a constant scalar curvature \K\
metric. The aim of the present paper is to extend these results to
the framework of extremal metrics.

\section{Statement of the result}

Let $(M,J, \omega)$ be a \K\ manifold with complex structure $J$ and
\K\ form $\omega$ and let $g$ denote the metric associated to the
\K\ form $\omega$, so that
\[
\omega (X,Y) = g (J \, X, Y).
\]
Further assume that $g$ is an extremal metric. Since the
automorphism group of any blow up of $M$ can be identified with a
subgroup of $\mbox{Aut}(M,J)$, and in light of the above mentioned
result of Calabi-LeBrun-Simanca about the isometry group of any
extremal metric, our strategy is to fix {\em a priori} a compact
subgroup $K$ of $\mbox{Isom}(M,g)$ and work $K$-equivariantly. Such
a $K$ will then be contained in the isometry group of the extremal
metric we are seeking on the blow up of $M$ at any set of points
$p_1, \ldots, p_n \in M$ in the $\mbox{Fix} \, (K_0)$, the fixed locus
of the identity component  $K_0$ of $K$. 
We will denote by $\frk$ the Lie algebra associated to
identity component of $K$. Observe that elements of $\frk$ vanish at
the points $p_1, \ldots, p_n$ to be blown up and hence these vector
fields can be lifted to the blown up manifold.

In order to produce extremal metrics on the blown up manifold, we
have to identify, among all $C^{\infty}$ functions on the blown up
manifold, those who generate real-holomorphic vector fields, since
these can arise as scalar curvatures of extremal metrics. 
To this aim, we define $\frh$ to be the vector space of
$K$-invariant hamiltonian real-holomorphic vector fields on $M$ 
or equivalently, the Lie algebra of the group $H$ of exact simplectomorphisms
commuting with $K$. The
correspondence between real-holomorphic vector fields and the scalar
functions on $M$ can be encoded in a compact way in a moment map
$\xi_\omega$
\[
\xi_\omega \colon M \rightarrow \frh^*
\]
for the action of $H$ uniquely determined by imposing to have mean
zero. More explicitly the function $f : = \ip{\xi_\omega , X}$
associated to the vector field $X \in \frh$ is defined to be the
unique solution of
\[
- \rd f  =  \omega (X, -)
\]
whose mean value over $M$ is $0$.

Using this setup, our main result reads~:
\begin{teor}
\label{mainthm-2} Let $(M,J,\omega)$ be a compact $m$-dimensional
\K\ manifold whose associated  \K\ metric $g$ is extremal, and let $K$ be
a compact subgroup of $\mbox{Isom}(M,g)$ whose Lie algebra contains
the vector field $J \, \nabla {\bf s}$ as well as any element of
$\frh$. Let $K_0$ denote the identity component of $K$.

Given $p_1, \dots ,p_n \in \mbox{Fix}\, (K_0)$ and $a_1, \ldots, a_n
>0$ such that $a_{j_1} = a_{j_2}$ if $p_{j_1}$ and $p_{j_2}$ are
in the same $K$-orbit, there exists $\e_0>0$ and, for all $\e \in (0, \e_0)$, there
exists a $K$-invariant extremal \K\ metric $\omega_\e$ on $\tilde M$,
the blow up of $M$ at $p_1, \ldots, p_n$, such that its associated \K\ form $\omega_\e$
lies   in the class
\[
 \pi^*[\omega] -  \e^{2} \, \left(a_1^{\tfrac{1}{m-1}} \, PD [E_1] +
\ldots + a_n^{\tfrac{1}{m-1}} \, PD \, [E_n] \right)
\]
where $\pi \colon \tilde M \rightarrow M$ is the standard projection
map, the $PD[E_j]$ are the Poincar\'e duals of the $(2m-2)$-homology
classes of the exceptional divisors of the blow up at $p_j$.

Finally, the sequence of metrics $(g_\e)_\e$ converges to $g$ (in
smooth topology) on compacts, away from the exceptional divisors.
\end{teor}

It is important to stress that our analytical construction does not give {\em one}
extremal metric but {\em a family} converging to the starting metric on the base manifold.
For such a construction to work, it is necessary to have  
$p_1, \dots ,p_n \in \mbox{Fix}\, (K_0)$ and
$J \, \nabla {\bf s} \in \frk$. On the other hand the condition $\frh \subset \frk$, while often
satisfied in important examples such as toric manifolds with $K$ giving the torus action, is certainly  far from being necessary.
We give a simple geometric sufficient condition for it to hold
in Proposition~\ref{fixx}.

In the general case when  $\frh$ is not included in $\frk$, there is a
natural decomposition
\[
\frh = \frh' \oplus \frh'' ,
\]
where $\frh' : =  \frh \cap \frk$ is the subspace of $K$-invariant
real-holomorphic vector fields in $\frk$. The previous result then
appears as a special case of the more general~:
\begin{teor}
\label{mainthm} Assume that $(M, J, \omega)$ is a compact \K\ manifold whose associated  \K\ metric $g$ is extremal, and let $K$ be a
compact subgroup of $\mbox{Isom}(M,g)$ whose Lie algebra contains
the vector field $J \, \nabla {\bf s}$. Let $K_0$ denote the identity component of $K$. 
We decompose the space of
$K$-invariant hamiltonian  real-holomorphic vector fields $\frh =
\frh' \oplus \frh''$ where $\frh' =(\frh \cap \frk)$. Given $p_1,
\dots ,p_n \in \mbox{Fix} \, (K_0)$ such that~:
\begin{itemize}
\item[(i)] there exists $a_1, \ldots, a_n  >0$ satisfying
\[
\sum_{j} a_j \, \xi_\omega (p_j) \in \frh' \,^* ,
\]
and $a_1, \ldots, a_n
>0$ such that $a_{j_1} = a_{j_2}$ if $p_{j_1}$ and $p_{j_2}$ are
in the same $K$-orbit,

\item[(ii)] the projections of  $\xi_\omega (p_1), \ldots,
\xi_\omega (p_n)$ over $ \frh''  \,^*$ span $\frh''  \,^*$,\\[3mm]

\item[(iii)] there is no nontrivial elements of $\frh''$ that vanishes at
$p_1, \ldots, p_n$,

\end{itemize} there exists $\e_0 > 0$ and, for all $\e \in
(0, \e_0)$, there exists a $K$-invariant extremal \K\ metric
$g_\e$ on $\tilde M$, the blow up of $M$ at $p_1, \ldots, p_n$,
whose associated \K\ form $\omega_\e$ lies in the class
\[
\pi^*[\omega] -  \e^{2} \, \left( a_1^{\tfrac{1}{m-1}} \, PD [E_1] +
\ldots + a_n^{\tfrac{1}{m-1}} \, PD \, [E_n] \right)
\]
where $\pi \colon \tilde M \rightarrow M$ is the standard projection
map, the $PD[E_j]$ are the Poincar\'e duals of the $(2m-2)$-homology
classes of the exceptional divisors of the blow up at $p_j$.

Finally, the sequence of metrics $(g_\e)_\e$ converges to $g$ (in
smooth topology) on compacts, away from the exceptional divisors.
\end{teor}

When $\frh \subset \frk$,  $\frh'= \frh$ and 
$\frh''=\{0\}$, hence (i), (ii) and (iii) become vacuous and
Theorem~\ref{mainthm} reduces to Theorem~\ref{mainthm-2}. 
In \S 4 we explain how conditions (i)-(iii) arise in our analytical approach.

\begin{rmk} Condition (iii) can be removed if we leave some freedom
on the weights of the exceptional divisor on the blown up manifold.
More precisely, Theorem~\ref{mainthm} still holds without assuming
(iii) but in this case, the only information we have about
$[\omega_\e]$ reads
\[
\omega_\e \in  \pi^*[\omega] -  \e^{2} \, \left( \tilde a_1^{\tfrac{1}{m-1}} \, PD
[E_1] + \ldots + \tilde a_n^{\tfrac{1}{m-1}} \, PD \, [E_n] \right)
\]
where $\tilde a_1, \ldots, \tilde a_n >0$ depend on $\e$ and satisfy
\[
|\tilde a_j - a_j |\leq c\, \e^{\tfrac{2}{2m+1}}.
\]
In other words, by removing (iii) we slightly loose the control on
the \K\ classes.
\end{rmk}

Theorem~\ref{mainthm} is a generalization of the constructions given
in \cite{ap} and \cite{ap2}. Indeed, in \cite{ap} is treated the
case where $g$ is a constant scalar curvature \K\ metric $K =
\{Id\}$ and where $\frh = \{0\}$, while in \cite{ap2} is treated the
case where $g$ is a constant scalar curvature \K\ metric, $K$ is a
discrete subgroup of $\mbox{Isom} \, (M,g)$, $\frh' = \{0\}$, and
$\frh''$ is not necessarily trivial.

The choice of the symmetry groups $K$ is a really delicate problem.
Indeed, given the fact that the blown up points have to be chosen in
$\mbox{Fix} \, (K_0)$ it is rather natural to choose $K_0$ to be fairly
small, so that its fixed-point set is large. However, the smaller
$K_0$ the larger $\frh$ and hence the harder it is to fulfill the
requirement that $\frh \subset \frk$ in Theorem~\ref{mainthm-2} or
the requirement that conditions (i) and (ii) in
Theorem~\ref{mainthm} are fulfilled. Conditions (i) and (ii) are of
course difficult to check given the fact that the moment map
$\xi_\omega$ is in general hard to analyze. Nevertheless, there are
large classes of manifolds, notably all toric ones, for which
computations can be done. Concrete examples listed in section \S 11
will illustrate how delicate this issue is.

Once the necessary notations and definitions are introduced, in \S 4 we will show how
conditions (i)-(iii) naturally arise in our construction of the converging family of extremal
metrics $g_{\e}$.

\begin{rmk} Condition (i) and (ii) should be related to Mabuchi's $T$-stability
\cite{ma} and to Szekelyhidi's relative $K$-stability \cite{sz} in the
same way the analogue conditions for constant scalar curvature
metrics are related to the asymptotic Chow semi-stability along the
line of ideas described by Thomas in \cite{th} (pages 27 and 28).
Indeed, if instead of fixing the group $K$ a priori, we fix the set of points $\{p_1,\dots, p_n\}$
to be blown up, a natural choice of $K$ would be any maximal torus in the subgroup
of ${\mbox{Isom}}_0(M,g)$ fixing each $p_j$. We believe that with these choices,
conditions (i) and (ii)  should be equivalent to the relative $K$-stability of the 
blown up manifold, when
the resulting \K\ classes are rational (which in turn should be equivalent to a relative GIT
stability of the configurations of points in $M^n$ \cite{dv}). 
Moreover let us observe that one can apply our construction to {\em any} extremal representative
of the class $[\omega]$. While loosing control on the explicit shape of the metric $g_{\e}$,
this clearly allows to find families of extremal representatives in $[\omega_{\e}]$, and gives
some flexibility on the choice of points and weights for which our construction works 
{\em for some} representative in $[\omega]$. This flexibility is indeed connected to
the above mentioned stability question, and it will be investigated in detail in \cite{aps}.
A first simple appearance of this freedom will be used below in the case of projective spaces.
\label{stopen}
\end{rmk}

If the initial manifold has constant scalar curvature, it might well
be that the extremal metrics we obtain are in fact constant scalar
curvature metrics. There is a simple criterion involving the points
$p_1, \ldots , p_n$ and the parameters $a_1, \ldots, a_n$, which
ensures that this is not the case.
\begin{prop}
Under the assumptions of Theorem~\ref{mainthm} (or
Theorem~\ref{mainthm-2}), if $ \sum_j a_j \, \xi_\omega (p_j) \neq
0$ then the metrics we obtain on $\tilde M$ are extremal with
nonconstant scalar curvature. \label{mainprop}
\end{prop}

We now emphasize the consequences of the above results for
projective spaces and more generally for toric varieties.

When
$(M,\omega)$ is ${\mathbb P}^{m}$ endowed with the \K\ form
$\omega_{FS}$ associated to a Fubini-Study metric, we let $(z^{1},
\ldots, z^{m+1})$ be complex coordinates in ${\mathbb C}^{m+1}$ and {\em let us fix for 
the rest of the paper the convention that $[\omega_{FS}] = PD[{\mathbb P}^{m-1}]$,
where ${\mathbb P}^{m-1}\subset {\mathbb P}^{m}$ is a linear subspace}. 
This is particularly relevant when getting quantitative estimates on the \K\ classes
reachable by our constructions. 
We
consider the group $K: = S^{1} \times \ldots \times S^{1}$, the
maximal compact subgroup of $PGL(m+1)$, whose action is given by
\[
 \begin{array}{cccclllll}
K  \times {\mathbb P}^m
 & \longrightarrow & {\mathbb P}^m \\[3mm]
\left( (\alpha_1, \ldots, \alpha_{m+1}), [z^{1}, \ldots , z^{m+1}]
\right) & \longmapsto & [\alpha_1 z^{1}, \ldots, \alpha_{m+1}
z^{m+1} ]
\end{array}
\]
and we consider the set of fixed points of $K$
\[
p_1 : =  [1: 0 : \ldots : 0], \quad \ldots  ,  \quad p_{m+1}: = [0 :
\ldots : 0 :  1]
\]
In this case, the space $\frh$ is spanned by vector fields of the
form
\[
\Re  \, \left( z^{j} \, \del_{z^j} - z^k \,\del_{z^k} \right)
\]
and we have $\frk = \frh = \frh'$ and $\frh'' = \{0\}$. As a
consequence of the result of Theorem~\ref{mainthm-2}, we obtain
extremal \K\ metrics on the blow up of ${\mathbb P}^m$ at the points
$p_1, \ldots, p_n$, for any $n=1, \ldots, m+1$.

It is worth emphasizing that the special structure of the points
which can be blown up on ${\mathbb P}^m$ has its origin in the fact
that we are starting from a specific choice of a Fubini-Study metric
and hence, away from the blow up points the extremal \K\ metric
$\omega_\e$ is close to $\omega_{FS}$. 
This example shows well the {\em riemannian} nature of our results
(see Remark ~\ref{stopen}).
Now, if $q_1, \ldots, q_n \in
{\mathbb P}^m$ are linearly independent one can find extremal
metrics on the blow up of ${\mathbb P}^m$ at $q_1, \ldots, q_n$ but
this time the metric will be close to $\psi^* \omega_{FS}$ away from
the blow up points, where $\psi$ is an automorphism of the
projective space such that
\[
\psi  (p_j) = q_j  .
\]
Yet, since $[\psi^* \omega_{FS}]$ is independent of $\psi$ and of
the choice of the Fubini-Study metric, we have obtained the following {\em \K ian}
version of Theorem ~\ref{mainthm-2} for ${\mathbb P}^m$~:
\begin{corol}
Fix $1 \leq n\leq m+1$. Given $q_1, \ldots, q_n \in {\mathbb P}^m$
linearly independent points and $a_1, \ldots, a_n
>0$, there exists $\e_0 >0$ and for all $\e \in (0, \e_0)$ there
exists an extremal \K\ metric $g _\e$ on the blow up of ${\mathbb
P}^m$ at $q_1, \ldots, q_n$ whose associated \K\ form $\omega_\e$ lies in the class
\[
\pi^*[\omega_{FS}] -  \e^{2} \, \left( a_1^{\tfrac{1}{m-1}} \, PD[E_1]
+ \ldots + a_n^{\tfrac{1}{m-1}} \, PD [E_n] \right)
\]
In addition, the \K\ metrics $g_\e$ do not have constant scalar
curvature unless $n = m+1$ and $a_1 =\ldots = a_{m+1}$.
\label{coproj}
\end{corol}

The conditions $n = m+1$ and $a_1 =\ldots = a_{m+1}$ being necessary and 
sufficient to get 
constant scalar curvature metrics among our family of extremal ones fits exactly with 
the more familiar picture of the Futaki invariants. Calabi has in fact proved that 
an extremal metric has constant scalar curvature iff its Futaki invariant vanishes \cite{ca2},
and we will show in \S 11, using Mabuchi's result \cite{ma1} relating the 
Futaki invariant to the 
coordinates of the barycenter of the convex polytope of a toric variety, that the above 
conditions are indeed equivalent to the vanishing of the
Futaki invariants for blow ups of  ${\mathbb P}^m$.

The case corresponding to $n=1$ in Corollary~\ref{coproj} was
already obtained by Calabi in more generality (i.e. for all \K\
classes) \cite{ca} and the case where ${\mathbb P}^m$ is blown up at
$m+1$ linearly independent points $q_1, \ldots, q_{m+1}$ and $a_1
=\ldots = a_{m+1}$ was already studied in \cite{ap2} where constant
scalar curvature metrics were obtained.

In the case where ${\mathbb P}^m$ is blown up at more than $m+1$
points in general position the resulting manifolds do not have
nonzero holomorphic vector fields, hence extremal metrics are forced
to have constant scalar curvature and the existence of some constant
scalar curvature \K\ metrics follows from \cite{ap2} and
\cite{Rol-Sin}.

The previous Corollary can be understood as a special case of the
the existence of extremal metrics on the blow up of toric varieties, which in fact 
leads to a more general result as we will see below even for ${\mathbb P}^m$.
If $(M,J, \omega )$ is a $m$-dimensional toric variety whose associated metric is extremal, one can take
$K$ to be the torus $T^m$ giving the torus action. It then follows from Proposition 
~\ref{fixx}
that $\frh =  \frk$, the Lie algebra associated to $K$. One can
apply Theorem~\ref{mainthm-2} to get~:
\begin{corol}
Assume that $(M, J, \omega)$ is a toric variety whose associated metric is extremal,
and let $K$ be the torus $T^m$ giving the torus
action. Given $p_1, \ldots, p_n \in \mbox{Fix} \, (K)$ and $a_1,
\ldots, a_n > 0$, there exists $\e_0 > 0$ and for all $\e \in (0,
\e_0)$ there exists an extremal \K\ metric $g_\e$ on the blow up
of $M$ at $p_1, \ldots, p_n$ whose associated \K\ form $\omega_\e$ lies in the class
\[
\pi^*[\omega] -  \e^{2} \, \left( a_1^{\tfrac{1}{m-1}} \, PD[E_1] +
\ldots + a_n^{\tfrac{1}{m-1}} \, PD [E_n] \right)
\]
\label{coproj-2}
\end{corol}
In other words, one can blow up {\em any set of points contained in
the fixed-point set of the torus-action} and the weights $a_j
>0$ can be chosen arbitrarily.

Since blowing up a toric variety at such points preserves the toric
structure, one can apply inductively the last Corollary. Therefore,
we obtain extremal metrics on any such iterated blow up. Beside
these applications this last Corollary can be applied to can be
applied to any toric \K-Einstein manifold, the classification of
which has been completed in dimension $m = 2, 3$ and $4$, in
\cite{ma}, \cite{nak1} and all the symmetric examples have been found
by Batyrev and Selivanova \cite{bs} and to the
one parameter family of extremal metrics found by Calabi of the blow
up of ${\mathbb P}^m$ at one point, producing then a wealth of open subsets
of classes in the \K\ cone which have extremal representatives.

For example our result applied to  ${\mathbb P}^{2}$, $ {\mathbb P}^{1}  
\times {\mathbb P}^{1}$ and $Bl_p{\mathbb P}^{2}$ as base 
manifolds leads to the following:

\begin{corol}
\begin{enumerate}
\item
If $M=Bl_{p_1,p_2}{\mathbb P}^{2}$ then the following \K\ classes have extremal representatives
\[
\pi^*[\omega_{FS}] -  \, \left( a_1 \, PD[E_1] + \e^{2}  a_2 \, PD [E_2] \right),
\qquad  a_1 < 1
\]

\[
\pi^*[\omega_{FS}] -  \, \tfrac{a_1 - \e^2}{a_1+a_2-\e^2}\, PD[E_1] -
\tfrac{a_2 - \e^2}{a_1+a_2-\e^2}\,  \, PD [E_2], 
\]

\item
If $M=Bl_{p_1,p_2,p_3}{\mathbb P}^{2}$ and the points do not lie on a complex line, then the following \K\ classes have extremal representatives
\[
\pi^*[\omega_{FS}] -  \, \left( a_1\, PD[E_1] + \e^{2}  a_2 \, PD [E_2] 
+ \e^{2}  a_3 \, PD [E_3] \right),
\qquad  a_1 < 1
\]

\[
\pi^*[\omega_{FS}] -  \, \tfrac{a_1 - \e^2}{a_1+a_2-\e^2}\, PD[E_1] -
\tfrac{a_2 - \e^2}{a_1+a_2-\e^2}\,  \, PD [E_2] - \e^4 a_3PD[E_3],
\]

\[
\pi^*[\omega_{FS}] -  \, \tfrac{1 -\e^2(a_1 +a_2)}{2-\e^2(a_1+a_2+a_3)}\, PD[E_1] -
\, \tfrac{1 -\e^2(a_1 +a_3)}{2-\e^2(a_1+a_2+a_3)}\, PD[E_2] -
\, \tfrac{1 -\e^2(a_2 +a_3)}{2-\e^2(a_1+a_2+a_3)}\, PD[E_3]  
\]
where $a_j +a_k < 1$ for all $j, k$, and $a_1+a_2+a_3 <2$.\\

\item
If $M=Bl_{p_1,p_2,p_3}{\mathbb P}^{2}$ and the points lie on a complex line, then the 
following \K\ classes have extremal representatives
\[
\pi^*[\omega_{FS}] -  \, \left( a_1\, PD[E_1] + \e^{2}  a_2 \, PD [E_2] 
+ \e^{4}  a_3 \, PD [E_3] \right),
\qquad  a_1 < 1
\]

\[
\pi^*[\omega_{FS}] -  \e^2 \, \left( a\, PD[E_1] +  a\, PD [E_2] 
+ b \, PD [E_3] \right),
\qquad  b  < a.
\]

\end{enumerate}
\label{claproj}
\end{corol}

\begin{rmk}
This last family of examples is interesting also because it has been shown 
by A. Della Vedova \cite{dv}, building on Szekelyhidi's work, that the above \K\ classes
 do not have extremal representatives for $b>2a$, giving then an explicit upper bound
 for our construction to work.
 \end{rmk}
 
In all the above cases, the first families of classes are immediately obtained by our direct construction applied once or twice to  $Bl_p{\mathbb P}^{2}$ with a Calabi's metric. The other classes are obtained by applying our result either to $ {\mathbb P}^{1}  
\times {\mathbb P}^{1}$ with a product of Fubini-Study metrics, or by using some classical algebraic constructions which will be recalled in \S 11.

We should also recall that when we blow up three not aligned points, the \K\ classes
$ \pi^*[\omega_{FS}] -  \, \left( a_1\, PD[E_1] + a_2 \, PD [E_2] 
+  a_3 \, PD [E_3] \right)$,
where all the $a_j$ are sufficiently close to $\tfrac{1}{3}$, also have extremal representatives, thanks to the existence of a \K -Einstein metric on the resulting manifold, as shown by Siu-Tian-Yau, and the recalled deformation theory of Lebrun and SImanca \cite{ls}.

\section{Notation and conventions}

The following conventions are used throughout. If $(M,J)$ is a
complex manifold, we write $\mbox{Aut}(M,J)$ for the group of
biholomorphic maps $M\to M$. If $(M,\omega)$ is a symplectic
manifold, we write $\mbox{Exact} (M,\omega)$ for the group of exact
symplectomorphisms; that is, those that are generated by hamiltonian
vector fields. Finally if $(M,g)$ is a riemannian manifold, we write
$\mbox{Isom} (M,g)$ for the group of isometries of $(M,g)$.  We
denote by a subscript $0$ the identity-component of these groups
(even though the group of exact symplectomorphisms is already
connected).

The metric $g$, K\"ahler form $\omega$ and complex structure $J$ are
related by
\begin{equation}\label{e1.25.11.5}
g(JX,Y) = \omega(X,Y),\qquad \mbox{or equivalently} \qquad
\omega(X,JY) = g(X,Y).
\end{equation}

The action of $J$ commutes with the musical isomorphisms, but if
$\alpha$ is a $1$-form and $X$ is a vector field, we have
\begin{equation}\label{e2.25.11.5}
   J \, \alpha(X) = - \alpha (J X).
\end{equation}
Then $T^{1,0}$ corresponds to the $+i$-eigenspace of $J$ while
$\Lambda^{1,0}$ corresponds to the $-i$-eigenspace of $J$. In
particular, we have
\begin{equation}
\bar \del  f = \tfrac{1}{2} (\rd f  - i J \rd f),\qquad  \qquad J\,
\rd f = i(\bar \del f - \del f) \label{e3.25.11.5}
\end{equation} and
so on.

Recall that a vector field $X$ is said to be a {\em hamiltonian
vector field} if there exists a smooth {\em real valued} function
$f$ satisfying
\begin{equation}
\label{e1.26.11.5} X = J \nabla f \, .
\end{equation}
In this case we will write $X =X_f$. Using (\ref{e1.25.11.5}) we see
that this equation is always equivalent to
\begin{equation}\label{e2.26.11.5}
\omega \, ( X_f , - ) = - \rd f.
\end{equation}
or, using (\ref{e3.25.11.5}), is also equivalent to
\begin{equation}
\label{eq:ldld}
 \tfrac{1}{2} \, \omega ( \Xi_f, -) = - \bar \del f.
\end{equation}
when~\footnote{To help the reader connecting this notation with the
existing literature, let us remark that $\Xi_f = \, 2 \, i \,
\bar{\partial}^{ \uparrow } f = 2 \, i  \, \bar{\partial}^{\#} f$,
the $(1,0)$ part of the gradient of $f$, in the notations used in
\cite{ca2} and \cite{ls}.}
\[
\Xi_f := X_f  - i J X_f \in T^{1,0}
\]

Let us now define the second order operator
\begin{equation}\label{e3.26.11.5}
\begin{array}{rccccllll}
P_\omega : & \ci(M) & \longrightarrow  & \Lambda^{0,1}(M,T^{1,0}),\\[3mm]
& f &  \longmapsto  &  \tfrac{1}{2} \, \bar \del \, \Xi_f
\end{array}
\end{equation}
so that the null-space of $P_\omega$ (beside the constant function)
corresponds to holomorphic vector fields with zeros. Observe that
the operator $P_\omega$ depends on the \K\ metric $\omega$. Also,
with this definition, a metric $\omega$ is extremal if and only if
$P_\omega({\bf s} (\omega)) = 0$.

Clearly, any smooth, complex valued function $f$ solution of
\[
P^*_\omega \, P_\omega \, f = 0
\] on $M$ gives rise to a holomorphic vector field $\Xi_f$ defined
by (\ref{eq:ldld}) since by integration over $M$ implies that
$\|\bar \del \, \Xi_f \|_{L^2(M)} =0$). We recall the following
important result which shows that the converse is also true~:
\begin{prop} \cite{ca2}, \cite{ls}
\label{compoten}  $\Xi \in T^{1,0}$ is a holomorphic vector field
with zeros if and only if there exists a {\em complex valued}
function $f$  solution of $P^*_\omega \, P_\omega  \, f =0$ such
that $\tfrac{1}{2} \, \omega (\Xi, -) = - \bar \del f$.
\end{prop}

In addition, we have the following result which follows from a
theorem of Lichnerowicz (see Besse \cite{be}, Corollary~2.125 and
\cite{ls})
\begin{prop}
\cite{be}, \cite{ls} A vector field $X$ is a Killing vector field
with zeros if and only if there exists a {\em real valued} function
$f$ solution of $P^*_\omega \, P_\omega \, f =0$ such that $\omega
(X, -) = - \rd f$. \label{killhamilton}
\end{prop}

In other words, if $\Xi$ is a holomorphic vector field and  $f$ the
function given in Proposition \ref{compoten}, then $f$ can be chosen
to be real valued when $X = \Re \,\Xi$ is a Killing vector field.
Also, any Killing vector field is automatically real-holomorphic.

Observe that, in particular, if $X$ is a Killing vector field with
zeros, then $\Xi = X -i J X$ is a holomorphic vector field. We
recall the~:
\begin{defin} A vector field $X$ is real-holomorphic if and only if $X - iJX$ is a {\em holomorphic}
section of $T^{1,0}M$.
\end{defin}

\section{Equivariant set-up}

We fix a compact subgroup $K$ of $\mbox{Isom} (M, g)$ and we assume
that the Lie algebra of $K$ contains $X_{\bf s} = J \, \nabla {\bf
s}$. We do not insist that $K$ be connected.

Let us denote by $\frh$ the Lie algebra of real-holomorphic vector
fields which are $K$-invariant and are hamiltonian. Note that, since
$\omega$ is $K$-invariant, $X_{\bf s}$ certainly lies in $\frh$ for
any choice of $K$.

There is a large flexibility in the choice of $K$ varying from the
two extreme cases where $X_{\bf s}$ happens to generate a closed
subgroup of $\mbox{Isom}_0 (M,g)$, then this will be a
circle-subgroup $S$ contained in the center of $\mbox{Isom}_0 (M,g)$
and one could choose $K=S$ and the opposite situation where we
choose $K = \mbox{Isom}_0 (M,g)$.

With slight abuse of notations, we will identify elements of $\frh$
with the real-holomorphic vector fields corresponding to the
infinitesimal action of $H$ on $M$. For any \K\ metric $\omega$,
denote by $\xi_ \omega$ the moment map for the action of $H$,
uniquely determined by requiring to have mean zero
\begin{equation}\label{e2.28.11.5}
\xi_\omega : M  \to  \frh^*
\end{equation}
Recall that this is defined as follows.  If $X \in \frh$, then the
function $f = \ip{\xi_\omega, X}$ on $M$ is a hamiltonian for the
vector field $X$, namely, a solution of~:
\begin{equation}
\omega(X, -) = - \rd \, f
\end{equation}
normalized by
\[
\int_M f \, \omega ^m = 0 .
\]
Observe that according to (\ref{e11.26.11.5}) we also have
\[
\tfrac{1}{2} \, \omega(\Xi , -) = - \bar \del \, \ip{\xi_\omega ,
X}
\]
where $\Xi = X - i \, J \, X$ is the holomorphic vector field
associated to $X$.

This is just an invariant way of introducing the potentials
corresponding to hamiltonian,  real-holomorphic vector fields with
zeros.

\begin{rmk}
Notice that as the K\"ahler form varies (among $K$-invariant forms
the moment map varies. In  section 4, we will explicitly study the
dependence of $\xi_\omega$ on $\omega$.
\end{rmk}

For the blow-up problem, we note that a vector field $X$ lifts to
$\tilde M$, the blow up of $M$ at $p_1, \ldots, p_n$, if and only if
it vanishes at each of the points $p_j$. If we have fixed the
isometry group to contain $K$, it follows that we only stand a
chance of blowing up points which are fixed by every element of $K$.
So, we suppose that
\begin{equation}
\mbox{For all  $j=1, \ldots, n$, $p_j \in \mbox{Fix} \, (K)$}.
\end{equation}

Now, if $\tilde{\omega}$ is a putative extremal K\"ahler metric on
$\tilde{M}$ its scalar curvature must be a sum of $K$-invariant
potentials corresponding to vector fields that vanish at the $p_j$
and are $K$-invariant and hence they have to correspond to vector
fields which are in $\frh'$. Thus we introduce the lie algebra
$\frh'$ that is given by
\[
\frh ' = \frk \cap \frh
\]

We denote by $\frh''$ the orthogonal complement of $\frh'$ in $\frh$
with respect to the scalar product
\[
(X, \tilde X)_{\mathfrak h} : = \int_M \ip{\xi_\omega , X} \,
\ip{\xi_\omega , \tilde X} \, dvol_g.
\]
Informally, potentials of the form $\ip{ \xi_\omega , X' }$ (for $X'
\in \frh'$) will the {\em good potentials} corresponding to vector
fields vanishing at the $p_j$, and hence lifting them to $\tilde{M}$
they can be used to deform the scalar curvature of the \K\ form
$\tilde{\omega}$. The potentials $\ip{\xi_\omega , X''}$ (for $X''
\in \frh''$) will be the {\em bad potentials} corresponding to
vector fields that do not necessarily lift to $\tilde{M}$ but, in
any case, these are potentials that will not be used in the
deformation of the scalar curvature of $\tilde{\omega}$.

To apply a perturbation argument, as in \cite{ap2}, we shall need to
solve two linear problems. First we need to find a function
$\Gamma$, a constant $\lambda$ and a vector field $Y' \in \frh'$
solutions of
\begin{equation}
 \tfrac{1}{2} \, P^*_\omega \, P_\omega \Gamma  + \ip{\xi_\omega ,
Y'} + \lambda = c_m \, \sum_{j=1}^n a_j \, \delta_{p_j}
\end{equation}
where the masses $a_j$ are positive and $c_m >0$ is a positive
constant only depending on the dimension $m$. The solvability of
this problem comes down to the {\em relative moment condition}~:
\begin{equation}
\sum_{j=1}^n a_j \, \xi_\omega (p_j) \in \frh' \,^* \mbox{ for some
}a_j>0
\end{equation}

Using this, we first consider a first perturbation of $\omega$, away
from the points to be blown up. This perturbed \K\ form we consider
is given explicitly by
\[
\hat \omega_\e : = \omega + i \, \del \, \bar \del ( \e^{2m-2} \,
\Gamma )
\]
where $\e >0 $ is a small parameter. This \K\ form is well defined
away from the points $p_j$ (provided $\e$ is chose small enough)
and, as will follow from the analysis in the next section, has
scalar curvature given by
\[
{\bf s} (\hat \omega_\e) = {\bf s} (\omega ) + \e^{2m-2} \, \left(
\ip{\xi_{\omega + i \, \del \, \bar \del ( \e^{2m-2} \, \Gamma )} ,
Y'} +  \lambda \right) + {\mathcal O} (\e^{4m-2})
\]

The final task will be to perturb this \K\ metric into an extremal
metric. To this aim,  given any (smooth) function $f$, we need to be
able to find a function $\phi$, a constant $\nu$, a vector field $Z'
\in \frh'$ and masses $b_j \in {\mathbb R}$ solutions of
\begin{equation}
\tfrac{1}{2} \, P^*_\omega \, P_\omega  \, \phi   + \nu +
\ip{\xi_\omega , Z'} +  c_m \, \sum_{j=1}^n b_j \, \delta_{p_j} = f
\end{equation}
The solvability of this problem is precisely equivalent to the {\em
genericity condition}~:
\begin{equation}
\mbox{The projections of }\xi_\omega (p_1), \ldots, \xi_\omega
(p_n)\mbox{ in }\frh''  \,^* \mbox{ span }\frh''  \,^*.
\end{equation}

The core of the paper is to show that these conditions are indeed
{\em sufficient} conditions to guarantee the existence of extremal
metrics in the appropriate classes.

\section{Linear operators}

\noindent The linearization of the mapping
\begin{equation}
\label{e5.25.11.5} f\longmapsto {\bf s} (\omega + i\del \bar \del f)
\end{equation}
is given by the formula
\begin{equation}\label{e6.25.11.5}
{\mathbb L} : =  - \tfrac{1}{2}  \, \Delta^{2}_g - \mbox{Ric}_g
\cdot \nabla^{2}_g
\end{equation}
where $\mbox{Ric}_g$ stands for the Ricci tensor of the metric $g$
associated to $\omega$. On the other hand,
\begin{equation}\label{e7.25.11.5}
P^*_\omega \, P_\omega = \Delta^{2}_g + 2 \, \mbox{Ric}_g \cdot
\nabla^{2}_g - J  \, X_{\bf s}  + i X_{\bf s} .
\end{equation}
where $P$ is the operator defined in \eqref{e3.26.11.5} and  $X_{\bf
s}$ is the hamiltonian vector field associated with ${\bf s}$.

Hence we have~:
\begin{equation} \label{e12.26.11.5} {\mathbb L} =
-\tfrac{1}{2} P^*_\omega \, P_\omega - \tfrac{1}{2} \, J  \,
X_{\bf s} + \tfrac{i}{2} \, X_{\bf s}
\end{equation}

Working equivariantly with respect to a compact group $K$ whose lie
algebra contains $X_{\bf s}$ has the important effect of  making the
last term in \eqref{e12.26.11.5} disappearing, leaving a real
operator on $K$-invariant functions.

Consider the map
\begin{equation}\label{e11.26.11.5}
\begin{array}{rcccc}
F : & \frh \times \ci(M)^{K} & \longrightarrow  & \ci(M)^{K}, \\[3mm]
  & (X, f) & \longmapsto & {\bf s} (\omega + i \del \bar \del f) - \ip{ \xi_{\omega
+i\del \bar \del f} , X }.
\end{array}
\end{equation}
Here the superscripts $K$ denote the $K$-invariant part of the
function space.

The following is due to Calabi and LeBrun--Simanca.
\begin{prop}
\label{linear} Assume that  $\omega$ is extremal and $ X_{\bf s} \in
\frh$, then $D_f F|_{(X_{\bf s},0)}$, the linearization  of $F$ with
respect to $f$ at $(X_{\bf s} ,0)$ is equal to $- \tfrac{1}{2}\,
P^*_\omega \, P_\omega$.
\end{prop}

\begin{proof}  We already know the linearization of the scalar
curvature map, so we only need to know the linearization of
\[
f \longmapsto \xi_{\omega + i \del \bar \del f}
\]
with respect to $f$. Take any $X \in \frh$. Since $f$ is
$K$-invariant, $X$ is a Killing vector field (with zeros) for the
\K\ form $\omega + i\del \bar \del f$. Hence, using the analysis of
\S 4, we can write
\[
\tfrac{1}{2} \, (\omega + i \del \bar \del  f)(\Xi , -) = - \bar
\del \ip{\xi_{\omega + i\del \bar \del f}, X}
\]
where $\Xi := X - i\, J \, X$, and we see immediately that
$\dot{\xi}$, the first variation of $f \longmapsto \xi_{\omega + i
\del \bar \del f}$ with respect to $f$ computed at $f=0$, satisfies
\[
\tfrac{i}{2} \, \del \bar \del \, f(\Xi, -) = - \bar \del \,
\ip{\dot{\xi} , X}.
\]
Working in local coordinates, the left hand side of this expression
is equal to
\[
\tfrac{i}{2}\, \rd \bar z^{j} \tfrac{\del}{\del \bar z^j}
\left(\Xi^k\tfrac{\partial f}{\partial z^k} \right)
\]
because $\Xi$ is holomorphic. Hence we see that
\begin{equation}\label{e5.26.11.5}
\ip{\dot{\xi} , X} = - \tfrac{i}{2}\, \Xi \, f.
\end{equation}
Now, we apply this analysis when $\omega$ is extremal, with extremal
vector field $X_{\bf s} \in \frh$. We obtain for any smooth function
$f$
\[
D_f F|_{(X_{\bf s}, 0)}( f ) = {\mathbb L} f + \tfrac{i}{2}\,
\Xi_{\bf s} \, f \qquad \qquad \mbox{ with } \qquad \qquad \Xi_{\bf
s} : = X_{\bf s} - i \, J \, X_{\bf s}.
\]
Hence
\begin{equation}
D_f F|_{(X_{\bf s}, 0)}( f ) =  -\tfrac{1}{2}  \, P^*_\omega \,
P_\omega f - \tfrac{1}{2} \, J \, X_{\bf s} \, f + \tfrac{i}{2} \,
X_{\bf s} \, f   + \tfrac{i}{2}\, \Xi_{\bf s}  \, f = -
\tfrac{1}{2} \, P^*_\omega \, P_\omega \,f + i \, X_{\bf s} \, f.
\end{equation}
Remembering that when $f$ is $K$-invariant and $X_{\bf s} \in \frk$,
we have
\[
X_{\bf s} f =0
\]
we conclude that $D_f F|_{(X_{\bf s}, 0)}( f ) = - \tfrac{1}{2} \,
P^*_\omega \, P_\omega \,f $. This completes the proof.
\end{proof}

\section{Burns-Simanca's metric on the blow up of ${\mathbb C}^m$ at the origin}

We describe a scalar flat \K\ form $\eta$ define on $\tilde{\mathbb
C}^m$, the blow up at the origin of ${\C}^m$. This metric is $U(m)$
invariant and was found by Burns \cite{leb}, when $m=2$, and Simanca \cite{sim}, when $m\geq 3$, following a method introduced in  \cite{ca} . 
Away from the exceptional divisor, the \K \ form $\eta$
is given by
\[
\eta  = i \, \del \, \bar \del E_m (v)
\]
where $v = (v^{1}, \ldots, v^m)$ are complex coordinates in
${\mathbb C}^m \setminus \{ 0 \}$ and where the function $E_m$ is
explicitly given, in dimension $m=2$, by
\[
E_2 (v) : =   \tfrac{1}{2} \, |v|^2 + \log  |v|^2
\]
while in dimension $m \geq 3$, even though there is no explicit
formula for $E_m$ we have the following expansion
\[
E_m (v) = \tfrac{1}{2} \, |v|^2 -  |v|^{4-2m} + {\mathcal O} (
|v|^{2-2m} )
\]
as $|v|$ tends to $\infty$. Observe that $E_m$ is defined up to a
constant. Details can be obtained either in \cite{ca}, or \cite{sim}
or even \cite{ap2}.

It is important to stress that the metric $\eta$ is defined in terms
of a choice of coordinates $(v^{1},\dots,v^{m})$, and any choice of
local coordinates around the point $p_j$ gives rise to a preferred
$\eta$. On the other hand the geometry of extremal metrics, and in
particular our choice of the group $K$, points to a preferred choice
as we will see in the next section.

\section{$K$-invariance and extensions on the blow up}

We discuss the crucial question of the lifting of objects (such as
the action of $K$, holomorphic vector fields, associated
potential,\ldots) all of which are defined on $M$, to the blown up
manifold.

Recall that, blowing up a $m$-dimensional complex manifold at a
point can be understood as a connected sum construction which can be
performed by excising a small ball in complex normal coordinates
around the point we want to blow up and replacing it by a large
neighborhood of the exceptional divisor in $\tilde \C^m$, the blow
up of $\C^m$ at the origin, keeping some compatibility between
metrics and complex structures on the different summands.

Now, $M$ is endowed with a $K$ invariant \K\ form $\omega$ and
$\tilde \C^m$ will be equipped with a suitable multiple of the
scalar flat \K\ form $\eta$ defined in the previous section. Since
we want the action of the group $K$ to lift to an isometric action
also for the new metrics we have to impose the condition that $K
\subset U(m)$ on the neighborhood of the point $p$ which will be
blown up, since $U(m)$ is the isometry group of any Burns-Simanca's metric
$\eta$. This is accomplished by linearizing on a small neighborhood
of $p$ the action of $K$, which is ensured by the following
classical result \cite{bm}~:
\begin{prop}
Let $D$ be a domain of a complex manifold and $G \subset \mbox{Aut}
\, (D, J)$ be a compact subgroup with a fixed point $p \in D$. In a
neighborhood of $p$, there exist complex coordinates centered at $p$
such that in these coordinates the action of $G$ is given by linear
transformations. \label{prop:link}
\end{prop}

We will refer to these  coordinates as {\em  $G$-linear
coordinates}. The following proposition, proved in \cite{ap2}, shows that one can find
$G$-linear coordinates which are also normal coordinates about $p$,
a fixed point of $G$.
\begin{prop}
Assume that $G \subset \mbox{Isom} \, (M,\omega)$ is compact. Then
there exist $(z^{1}, \ldots , z^{m})$, $G$-linear coordinates
centered at $p \in \mbox{Fix} \, (G)$ such that
\[
\omega  =  i \, \del \bar \del ( \tfrac{1}{2} \, |z|^2 + \varphi )
.
\]
where the function $\varphi$ is $G$ invariant and $\varphi =
{\mathcal O}(|z|^{4})$. \label{masterpiece}
\end{prop}

In our construction of extremal \K\ metrics on blow ups, this
proposition will be used in the following way. We apply the previous
result to $G=K$ close to a fixed point $p \in \mbox{Fix} \, (K)$. We
obtain normal coordinates, in $D$ a neighborhood of $p$ for which
the action of $K$ is linear. Given $X \in \frk$ a Killing vector
field vanishing at a point $p$ we can lift this vector field as a
vector field $\tilde X$ on $\tilde D$ the blow up of $D$ at the
point $p$. If $D$ is endowed with $\alpha \, \eta$, a multiple of
Burns-Simanca's metric, then $\tilde X$ will still be a Killing vector
field of $(\tilde D, J, \alpha \, \eta)$. In addition, $\tilde X$
still vanishes at some point on the exceptional divisor over $p$ as
is shown in the following~:
\begin{prop}
\label{vanlift} Let $X$ be a real-holomorphic vector field on $M$,
and let $p$ be any point in $M$ such that $X(p)=0$. We denote by
$\tilde{X}$ the lift of $X$ to the blow up of $M$ at $p$. Then,
there exists a point $q$ on the exceptional divisor over $p$ such
that $\tilde{X}(q) = 0$.
\end{prop}
\begin{proof}
For simplicity, we give the proof in the case where $m=2$. Given $z
: = (z^{1},z^{2})$ complex coordinates centered at $p$, we write
\[
X -iJX = X^{1} \, \partial_{z^{1}} + X^{2} \, \partial_{z^{2}},
\]
with $X^{i}(0)=0$. Let $(u^{1}, u^{2})$ be complex coordinates on
$\tilde \C^{2}$ such that $z^{1} = u^{1} \, u^{2}$ and $z^{2} =
u^{2}$, covering an affine chart of the exceptional divisor over
$p$. Then, $\tilde{X} - iJ\tilde{X}$, the lift of the vector field
$X-iJX$, is given by
\[
\tilde{X} - iJ\tilde{X} = \tfrac{z^{2} \, X^{1} - z^{1}
X^{2}}{(z^{2})^2} \,
\partial_{u^{1}} + X^{2}\partial_{u^{2}}.
\]
We can always write
\[
X^{1} = a \, z^{1} + b \, z^{2} + {\mathcal O} (|z|^2), \qquad
\mbox{and} \qquad X^{2} = c \, z^{1} + d \, z^{2} + {\mathcal O}
(|z|^2)
\]
for $z$ close to $0$.

We consider the point of the exceptional divisor corresponding to
the line $z^{1} = \lambda \, z^{2}$ (i.e. $(u^{1},
u^{2})=(\lambda,0)$). Obviously, we have
\[
\lim_{z^{2}\rightarrow 0} X^{2}(\lambda \, z^{2}, z^{2}) = 0
\]
for all $\lambda {\mathbb C}$ and
\[
\lim_{z^{2} \rightarrow 0} \tfrac{z^{2} \, X^{1} - z^{1}
X^{2}}{(z^{2})^2} \,(\lambda \, z^{2}, z^{2}) = - c\lambda^{2} +
(a-d)\lambda + b ,
 \]
Unless $c=0, d=a,b\neq 0$, in which case the point at infinity of
$\C =\{\lambda\}$ annihilates $\tilde{X}$, the equation $ -
c\lambda^{2} + (a-d)\lambda + b  =0$ has always a root $\lambda$
which corresponds to a zero of $\tilde X$.
\end{proof}
 
Let us end this section with a simple but useful result. Since the obstruction
for our construction to be successful is essentially contained in $\frh''$, it is interesting to
have some geometric efficient property which implies that $\frh''=0$.

To state this condition precisely, let us fix $p\in M$ and denote by 
$\rho \colon K \rightarrow Gl(T_pM)$
the representation of $K$ induced on $T_pM$ by the action of $K$ on $M$.

\begin{prop}
If there exists $p\in Fix(K)$ s.t. the maximal torus $T_{\frk}$ of 
$\rho(K)$ has dimension equal to $dim_{\C}M$, then $\frh''=0$.
\label{fixx}
\end{prop}
\begin{proof}
We need to show that $\frh \subset \frk$. If $p\in Fix(K)$ we can use Proposition
 linearize the action of 
$T_{\frk}^{\C}$ near $p$ in a such a way that $\rho(K)$ can be written as
$z_j \mapsto e^{i\theta_j} z_j$ in suitable complex coordinates on $T_pM$.
The condition $dim_{\R}T_{\frk} = dim_{\C}M$ then implies that $\theta_j \neq \theta_l$ for $j\neq l$. This immediately implies that the elements of $H$ are also diagonal, and hence
the result follows.
\end{proof}

Observe that the above condition is easily satisfied by any toric manifold with $K$ the 
maximal compact torus giving the torus action.

\section{Mapping properties}

For all $r > 0$, we agree that
\[
B_r : = \{ z \in {\mathbb C}^m \quad : \quad  |z| < r \},
\]
denotes the open ball of radius $r >0$ in ${\mathbb C}^m$, $\bar
B_r$ denotes the corresponding closed ball and
\[
\bar B_r^* : = \bar B_r -\{0\}
\]
the punctured closed ball. We will also define
\[
C_r : = {\mathbb C}^m - \bar B_r  \qquad \mbox{and} \qquad \bar C_r
: = {\mathbb C}^m - B_r
\]
to be respectively the complement in ${\mathbb C}^m$ of the closed
the ball and the open ball of radius $r >0$.

\subsection{Operators defined on $M - \{p_1, \ldots, p_n\}$}

Assume that we are given $n$ distinct points $p_1, \ldots, p_n \in
M$. For each $j = 1, \ldots, n$, we can choose complex coordinates
$z : = (z^{1}, \ldots, z^m)$ in a neighborhood of $0$ in ${\mathbb
C}^m$, to parameterize a geodesic ball of radius $r$ centered at
$p_{j}$ in $M$. Furthermore, as explained in the previous section,
these coordinates can be chosen to be normal at $p_j$ and to be
$K$-linear. In order to distinguish between the different
neighborhoods and coordinate systems, we agree that, for all $r$
small enough, say $r \in (0, r_0)$, $B_{j , r}$ (resp. $\bar B_{j,
r}$ and $\bar B_{j , r}^*$) denotes the open ball (resp. the closed
and closed punctured ball) of radius $r$ in the coordinates $z$
parameterizing a fixed neighborhood of $p_j$.

We fix $r_0$ small enough so that $\bar B_{j,r}$ are disjoint for
all $r \leq 4 \, r_0$. We set
\[
\bar M_{r} : = M - \cup_{j} B_{j,r}
\]
The weighted space for functions defined on the noncompact manifold
\begin{equation}
M^* : = M - \{ p_{1}, \ldots , p_{n} \}. \label{eq:MP0}
\end{equation}
is then defined as the set of functions whose decay or blow up near
any $p_{j}$ is controlled by a power of the distance to $p_{j}$.
More precisely, we have the~:
\begin{dfn}
Given $\ell \in {\mathbb N}$, $\alpha \in (0,1)$ and $\delta \in
{\mathbb R}$, we define the weighted space ${\mathcal C}^{\ell ,
\alpha}_{\delta} (M^*)$ to be the space of functions $f \in
{{\mathcal C}}^{\ell , \alpha}_{loc} (M^*)$ for which the following
norm is finite
\[
\| f \|_{{{\mathcal C}}^{\ell, \alpha}_{\delta} (M^*)}  : = \| f
\|_{{{\mathcal C}}^{\ell, \alpha} (\bar M_{r_0})} + \sum_{j=1}^n \,
\left( \sup_{r\leq r_0}  \, \left( r^{-\delta} \, \| f \,
_{|_{B_{j,r_0}}} (r\, \cdot \, ) \|_{{{\mathcal C}}^{\ell, \alpha}
(\bar B_2 - B_1 )} \right) \right).
\]
\label{de:MP1}
\end{dfn}

We are interested in the mapping properties of the operator
\[
L_\omega  : = - \tfrac{1}{2} P^*_\omega \, P_\omega
\]
which has been defined in \S 6. We define in $\bar B^*_{j,r_0}$ the
function $G_j$ by
\[
G_j (z) =  - \log \, |z|^2 \qquad \mbox{ when $m=2$ \qquad and}
\qquad G_j (z) = |z|^{4-2m} \qquad \mbox{ when $m\geq 3$}.
\]
Observe that, unless the metric $\omega$ is the Euclidean metric,
these functions are not solutions of the homogeneous equation
associated to $L_\omega$, however they can be perturbed into $\tilde
G_j$ solutions of the homogeneous problem $L_\omega \, \tilde G_j
=0$. Indeed, reducing $r_0$ if this is necessary, we know from
\cite{ap2} that there exist functions $\tilde G_j$ which are
solutions of $L_\omega \, \tilde G_j =0$ in $B_{j,r_0}^*$ and which
are asymptotic to $G_j$ in the sense that $\tilde G_j- G_j \in
{\mathcal C}^{4, \alpha}_{6-2m} (\bar B^*_{j,r_0})$ when $m \geq 4$
and $\tilde G_j - G_j \in {\mathcal C}^{4, \alpha}_{\delta} (\bar
B^*_{j,r_0})$ for any $\delta < 6-2m$ when $ m = 2 ,3$. The rational
behind these constructions is that the \K\ metric $\omega$ osculates
the Euclidean metric to order $2$ and hence
\[
L_\omega \, G_j \in {\mathcal C}^{0, \alpha}_{2-2m} (\bar
B_{j,r_0}^*).
\]

With the functions $\tilde G_j$ at hand, we define the {\em
deficiency spaces}
\[
{\mathcal D}_0 : = \mbox{Span} \{ \chi_1  ,  \ldots, \chi_n  \} ,
\qquad \mbox{and} \qquad {\mathcal D}_1 : = \mbox{Span} \{ \chi_1 \,
\tilde G_1 , \ldots, \chi_n \, \tilde G_n \} ,
\]
where $\chi_j$ is a cutoff function which is identically equal to
$1$ in $B_{j,r_0/2}$ and identically equal to $0$ in $M-B_{j,r_0 }
$. Furthermore, we assume that $\chi_j$ is $K$-invariant.

When $m \geq 3$, we fix $\delta \in (4-2m, 0)$ and when $m = 2$ we
choose $\delta \in (0,1)$. We define the operator
\[
\begin{array}{rclcllll}
{\mathcal L}_\delta : & ( {{\mathcal C}}^{4, \alpha}_{\delta}
(M^*)^K  \oplus {\mathcal D} ) \times {\mathfrak h}' \times {\mathbb
R}
& \longrightarrow & {{\mathcal C}}^{0,\alpha}_{\delta - 4} (M^*)^K \\[3mm]
& (f, X' , \mu )  & \longmapsto &  - \tfrac{1}{2} \, P^*_\omega \,
P_\omega f - \ip{\xi_\omega, X'} -  \mu ,
\end{array}
\]
Where ${\mathcal D} = {\mathcal D}_1$ when $m \geq 3$ and ${\mathcal
D} = {\mathcal D}_0 \oplus {\mathcal D}_1$ when $m=2$. The main
result of this section reads~:
\begin{prop}
Assume that the points $p_1, \ldots, p_n \in M$ are chosen so that $
$. Then, the operator ${\mathcal L}_\delta$ defined above is
surjective and has a kernel of dimension $ n + 1 + \mbox{dim} \,
{\mathfrak h}'$. \label{pr:f-5.1}
\end{prop}
\begin{proof} The proof of this result follows from the general
theorem described in \cite{Mel} and \cite{Maz}, nevertheless, we
choose here to describe an (almost) self contained proof. Recall
that, when working equivariantly with respect to a group $K$, the
kernel of $L$ is spanned by the functions of the form $\ip{\xi, X}
+\mu$ where $X \in {\mathfrak h}$ and $\mu \in {\mathbb R}$. Also
recall that, by assumption $\ip{\xi, X}$ has mean $0$.

Observe that ${\mathcal C}^{0, \alpha}_{\delta-4} (M^*)^K \subset
L^{1} (M)$ precisely when $\delta > 4-2m$. We use the fact that $L$
is self adjoint and hence, for $h \in L^{1} (M)$, the problem
\[
\tfrac{1}{2} \, P^*  P \, f  + \ip{\xi , X'} + \mu + \sum_{j=1}^n
b_j \, \delta_{p_j} = h
\]
is solvable in $W^{3,q}(M)$ for all $q \in [1, \tfrac{2m}{2m-1})$  if
and only if
\[
\mu \, \mbox{Vol}_g (M) + \sum_{j=1}^n b_j = \int_M h \, dvol_g
\]
(remember the $\ip{\xi, X}$ is normalized to have mean $0$) and, for
all $X'' \in \frh''$
\[
\sum_{j=1}^n b_j \, \ip{ \xi(p_j), X''}  = \int_M h \, \ip{\xi, X''}
\, dvol_g
\]
(remember that $\frh'$ and $\frh''$  are constructed so that
\[
(X', X'')_{\frh} : = \int_M \ip{\xi, X'} \,  \ip{\xi, X''} \, dvol_M
=0 ,
\]
for all $X' \in \frh'$ and all $X'' \in \frh''$). The first equation
gives the value of $\mu$ in terms of the $b_j$'s and the function
$h$. While, the second system in $b_j$'s is solvable since we have
assumed that the projection of $\xi(p_1), \ldots , \xi(p_n)$ over
$\frh'' \, ^*$ spans $\frh'' \, ^*$.

To complete the proof, we simply invoke regularity theory which
implies that $f \in {{\mathcal C}}^{4, \alpha}_{\delta} (M^*)^K
\oplus {\mathcal D}$. The estimate of the dimension of the kernel is
left to the reader since it will not be used in the paper.
\end{proof}

\subsection{Operators defined on $\tilde {\mathbb C}^m$}

We choose coordinates $u : =  (u^{1}, \ldots, u^m)$ to parameterize
$\tilde {\mathbb C}^m$, the blow up of ${\mathbb C}^m$ at the
origin, away from the exceptional divisor.

We start with the~:
\begin{dfn} Given $\ell \in {\mathbb N}$, $\alpha
\in (0,1)$ and $\delta \in {\mathbb R}$, we define the weighted
space ${{\mathcal C}}^{\ell , \alpha}_{\delta} (\tilde {\mathbb
C}^m)$ to be the space of functions $w \in {{\mathcal C}}^{\ell ,
\alpha}_{loc} (\tilde {\mathbb C}^m)$ for which the following norm
is finite
\[
\| f \|_{{{\mathcal C}}^{\ell , \alpha}_{\delta} (\tilde {\mathbb
C}^m)}  : = \| f \|_{{{\mathcal C}}^{\ell , \alpha} (\tilde {\mathbb
C}^m - C_1)} +  \sup_{r\geq R_0} \, r^{-\delta} \, \| f_{|_{C_2}} (
r \, \cdot \, ) \|_{{{\mathcal C}}^{\ell , \alpha} (\bar B_2 - B_1)}
.
\]
\label{de:MP2}
\end{dfn}

Given $\delta \in {\mathbb R}$, we define the operator
\[
\begin{array}{rclcllll}
\tilde {\mathcal L}_\delta : & {{\mathcal C}}^{4, \alpha}_{\delta}
(\tilde {\mathbb C}^m)^K & \longrightarrow &  {{\mathcal C}}^{0,
\alpha}_{\delta - 4} (\tilde {\mathbb C}^m)^K  \\[3mm]
& f & \longmapsto &  - \tfrac{1}{2} P^*_\eta  P_\eta \, f
\end{array}
\]
and recall the following result which is borrowed from \cite{ap2}
\begin{prop}
Assume that $\delta \in (0,1)$. Then the operator $\tilde {\mathcal
L}_\delta$ defined above is surjective and has a one dimensional
kernel spanned by a constant function. \label{pr:MP3}
\end{prop}

\subsection{Bi-harmonic extensions}

The following results are concerned with Bi-harmonic extensions
either on the complement of the unit ball or on the unit ball of
${\mathbb C}^m$ of boundary data defined on the unit sphere. Here
$\Delta$ denotes the Laplacian on ${\mathbb C }^m$.
\begin{prop}
Given $h \in {\mathcal C}^{4, \alpha} (\del B_1 )$, $k \in {\mathcal
C}^{2, \alpha} (\del B_1 )$ such that
\[
\int_{\partial B_1}  (4 \, m \, h - k) = 0
\]
there exists a function $W^{i} ( = W^{i}_{h,k} )\in {\mathcal C}^{4,
\alpha}_1 (\bar B_1^* )$ such that
\[
\Delta^{2} \,  W^{i} = 0 \qquad \mbox{in} \qquad B_1 \qquad \qquad
W^{i}= h \quad \mbox{and} \quad \Delta W^{i} = k \qquad \mbox{on}
\qquad \del B_1  .
\]
Moreover,
\[
\| W^{i} \|_{{\mathcal C}^{4, \alpha}_1 (\bar B_1^*)} \leq c \,
(\|h\|_{{\mathcal C}^{4, \alpha}(\del B_1 )} + \|k \|_{{\mathcal
C}^{2, \alpha}(\del B_1)})
\]

Given $h \in {\mathcal C}^{4, \alpha} (\del B_1 )$, $k \in {\mathcal
C}^{2, \alpha} (\del B_1 )$ such that
\[
\int_{\del B_1 } k =0
\]
there exists a function $W^{o} ( = W^{o}_{h,k} ) \in {\mathcal C}^{4
, \alpha}_{3-2m} ({\mathbb C}^m  - B_1)$ such that
\[
\Delta^{2} \, W^{o} = 0, \qquad \mbox{in} \qquad {\mathbb C}^m  -B_1
\qquad \qquad  W^{o} = h \quad \mbox{and} \quad \Delta W^{o} = k
\qquad \mbox{on} \qquad \del B_1 .
\]
Moreover,
\[
\| W^{o} \|_{{\mathcal C}^{4, \alpha}_{3-2m} (C_1 ) } \leq c \,
(\|h\|_{{\mathcal C}^{4, \alpha}(\del B_1 )} + \|k \|_{{\mathcal
C}^{2, \alpha}(\del B_1)})
\]
\label{pr:f-5.5}
\end{prop}
Let us briefly comment on the assumption. To this aim, let us
concentrate on the case where both $h$ and $k$ are constant
functions in which case their bi-harmonic extensions $W^{i}$ and
$W^{o}$ are given explicitly by
\[
W^{i} (z) = a  + b \, |z|^2
\]
and
\[
W^{o} (z)= c \, |z|^{4-2m} + d \, |z|^{2-2m}
\]
when $m \geq 3$ and
\[
W^{o} (z) = c \, \log |z| + d \, |z|^{-2}
\]
when $m=2$. It is easy to see that $W^{i} \in {\mathcal C}^{4,
\alpha}_1 (\bar B_1^*)$ if and only if $a=0$. While $W^{o}  \in
{\mathcal C}^{2, \alpha}_{3-2m} ({\mathbb C}^m  -B_1)$ if and only
if $ c=0$. These conditions lead to the constraints on the function
$h$ and $k$ as in the statement of the result.

\section{Nonlinear perturbation results}

\subsection{Perturbation of $\omega$}

The first perturbation we will perform is concerned with the
perturbation of the extremal \K\ form $\omega$ which is defined on
the manifold $M$. We keep the notations which have been introduced
in \S 6.

The \K\ metric $\omega$ being extremal, we have
\begin{equation}
{\bf s} (\omega)  =  \ip{\xi_{\omega} , X_{\bf s} } + \mu_{\bf s}
\label{eq:solutioninitial}
\end{equation}
for some Killing field
\[
X_{{\bf s}}\in {\mathfrak h}'
\]
and some constant $\mu_{\bf s} \in {\mathbb R}$.

Given $n$ points $p_1, \ldots, p_n \in M$ and real parameters $a_1,
\ldots, a_n > 0$. The points $p_j$ are precisely the points where
the manifold $M$ will be blown up and the parameters $a_j >0$ will
be closely related to the \K\ classes and also the volume of the
exceptional divisor on the blow up manifold, once the construction
will be complete.

We recall here the crucial assumption on the choice of the points
$p_j$ and the parameters $a_j$, namely, we assume that
\begin{equation}
\sum a_j \, \xi (p_j) \in {\mathfrak h}' \, ^*
\label{eq:condition-1}
\end{equation}
This condition is precisely the one which allows one to find a
function $\Gamma$, a vector field $Y' \in {\mathfrak h}'$ and a
constant $\lambda \in {\mathbb R} $  such that
\begin{equation}
\tfrac{1}{2} \, P^*_\omega P_\omega \, \Gamma + \ip{\xi_\omega, Y'}
+ \lambda = c_m \, \sum_{j=1}^n a_{j} \, \delta_{p_{j}}
\label{eq:gamma}
\end{equation}
where the constant $c_m$ is defined by
\[
c_m := 4 \, (m-1) \, (m-2) \, |S^{2m-1}| \qquad \mbox{ when $m \geq
3$ \quad and} \qquad c_2 : = 2 \, |S^{3}|
\]
Indeed, the existence of $\Gamma$ depends on the ability to choose
$Y' \in {\mathfrak h}'$ and $\lambda \in {\mathbb R}$ so that
\[
\ip{\xi_\omega , Y' } +  \lambda  -  c_m \, \sum_{j=1}^n a_{j} \,
\delta_{p_{j}}
\]
is "orthogonal" to the kernel of $- \tfrac{1}{2} \,P^*_\omega \,
P_\omega$. Since this kernel is precisely spanned by the functions
of the form $\ip{\xi, X} + \alpha$ where $X \in {\mathfrak h}$ and
$\alpha \in {\mathbb R}$, we see that (\ref{eq:condition-1}) is a
necessary and sufficient condition for the existence of $\Gamma$. We
also have the relation
\begin{equation}
\lambda =  c_m \, \sum_{j=1}^n  a_{j} \label{eq:lambda}
\end{equation}
and we know that the Killing field $Y' \in {\mathfrak h}'$ has to be
chosen so that
\begin{equation}
\int_M \ip{\xi_\omega , Y'} \,  \xi_\omega  \, dvol_g - \sum_{j=1}^n
a_j \, \xi_\omega (p_{j}) \in {\mathfrak h}'' \, ^*
\label{eq:Xprime}
\end{equation}
It is not hard to check that the function $\Gamma$ has a nice
expansion near each $p_{j}$. Indeed, we have the~:
\begin{lemma}
Near $p_{j'}$, the following expansions hold
\[
\Gamma (z) =  -  a_{j} \, |z|^{4-2m}  + {\mathcal O} (|z|^{6-2m})
\]
when $m\geq 4$,
\[
\Gamma (z) =  - a_{j} \, |z|^{-2} + b_{j} \  \log |z| + c_{j} +
{\mathcal O} (|z|)
\]
for some $b_{j} , c_{j} \in {\mathbb R}$ when $m =3$, while
\[
\Gamma (z) = a_{j} \, \log |z| + b_{j} + c_{j} \cdot z + {\mathcal
O} (|z|^2 \, (-\log |z|))
\]
for some $b_{j} \in {\mathbb R}$ and $c_{j} \in {\mathbb C}^m$, when
$m=2$. Here ${\mathcal O} (|z|^q  (-\log |z|)^{q'})$ denotes a
smooth function defined away from $0$, whose partial derivatives,
when taken with respect to the vector fields $|z|\, \partial_{z^j}$
and $|z|\, \partial_{\bar z^j}$, are bounded by a constant
(depending on the number of derivatives) times $|z|^q  (-\log
|z|)^{q'}$. \label{le:NPR1}
\end{lemma}

We fix
\begin{equation}
r_\e: = \e^{\tfrac{2m-1}{2m+1}} \label{eq:NPR1}
\end{equation}
This corresponds to the radius of the balls centered at the points
$p_{j}$ which will be excised from $M$. Recall that, for $r >0$
small enough we have defined
\[
\bar M_r : =  M - \cup_{j} B_{j,r}
\]
On each of the boundaries of $\bar M_{r_\e}$, we will need some
small boundary data. Hence, we assume that we are given \[ h_{j} \in
{\mathcal C}^{4, \alpha}(\del B_1)^K\quad \mbox{and} \qquad k_{j}
\in {\mathcal C}^{2, \alpha}(\del B_1)^K \] for  $j = 1, \ldots, n$,
satisfying
\begin{equation}
\| h_{j} \|_{{\mathcal C}^{4, \alpha} (\del B_1)} + \|
k_{j}\|_{{\mathcal C}^{2, \alpha} (\del B_1)} \leq \kappa \,
r_\e^{4} ,\label{eq:NPR2}
\end{equation}
where $\kappa >0$ will be fixed later on. The subscript $K$ in the
definition of the spaces is meant to remind the reader that the
functions $h_j$ and $k_j$ are invariant under the action of $K$. We
further assume that
\begin{equation}
\int_{\del B_1} k_{j} = 0 \label{eq:NPR3}
\end{equation}
so that the second half of the result of Proposition~\ref{pr:f-5.5}
applies. To keep notation short we set
\[
{\bf h} : = (h_1, \ldots, h_n) \qquad \mbox{and} \qquad {\bf k} : =
(k_1, \ldots, k_n).
\]

In particular, we can define the function $W_{\e, {\bf h}, {\bf k}}$
which is identically equal to $0$ in $\bar M_{2r_0}$ and, for $j=1,
\ldots, n$ is equal to
\begin{equation}
 W_{\e, {\bf h},
{\bf k}} : =   \chi_{j} \, W^{o}_{h_{j}, k_{j}} ( \cdot / r_\e) ,
\label{eq:doublev}
\end{equation}
in $B_{j,2r_0}$. Here $\chi_{j}$ is a cutoff function which is
identically equal to $1$ in $B_{j,r_0}$ and identically equal to $0$
in $M - B_{j,2r_0}$. Observe that the function $W_{\e, {\bf h}, {\bf
k}}$ depends (linearly) on ${\bf h}$ and ${\bf k}$ and is $K$
invariant.

This being understood, we have the~:
\begin{prop}
Given $\delta \in (4-2m, 5-2m)$ when $m\geq 3$ or $\delta \in (0,
2/3)$ when $m =2$, there exists $\e_\kappa >0$ and $c_\kappa
>0$ such that for all $\e \in (0, \e_\kappa)$ one can find a
function $\phi_{\e, {\bf h}, {\bf k}} \in {\mathcal C}^{4,
\alpha}(\bar M_{r_\e})$, a vector field $Y'_{\e, {\bf h}, {\bf
k}}\in \frh'$ and a constant $\lambda_{\e, {\bf h}, {\bf k}} \in
{\mathbb R}$ such that the scalar curvature of the \K\ form
\[
\omega_{\e, {\bf h}, {\bf k}} =  \omega + i  \, \del  \, \bar \del
\, ( \e^{2m-2} \, \Gamma + W_{\e, {\bf h}, {\bf k}} + \phi_{\e, {\bf
h}, {\bf k}})
\]
defined on $\bar M_{r_\e}$, satisfies
\[
-d {\bf s} (\omega_{\e, {\bf h}, {\bf k}}) =  \omega_{\e, {\bf h},
{\bf k}} (X_{\bf s} + \e^{2m-2} \, Y' + Y'_{\e, {\bf h}, {\bf k}} ,
-)
\]
with
\[
\tfrac{1}{{|\del B_{j, r_\e}|}} \int_{\del B_{j, r_\e}} {\bf s}
(\omega_{\e, {\bf h}, {\bf k}}) - {\bf s} (\omega_{\e, {\bf h}, {\bf
k}}) = \lambda_{\e, {\bf h} , {\bf k}}
\]
In addition, we have the estimate
\[
\begin{array}{rlllll}
\| Y'_{\e, {\bf h}, {\bf k}} \|_{L^\infty} + r_\e^{2m-4} \,
\sup_{j=1, \ldots, n} \, \| \phi_{\e, {\bf h},{\bf k}} \, |_{\bar
B_{j ,2 r_\e} -B_{j, r_\e}} (r_\e \,
\cdot) \|_{{\mathcal C}^{4, \alpha} ( \bar B_2 -B_1)}  \qquad \\[3mm]
\hfill \leq c_\kappa \, (r_\e^{2m+1} + \e^{4m-4} \, r_\e^{6 -4m
-\delta}),
\end{array}
\]
\[
|\lambda_{\e, {\bf h}, {\bf k}} | \leq c \, \e^{2m-2}
\]
And we also have
\begin{equation}
\begin{array}{llll}
|\lambda_{\e, {\bf h}^{(1)},{\bf k}^{(1)}} - \lambda_{\e, {\bf
h}^{(2)},{\bf k}^{(2)}} | + \| Y'_{\e, {\bf h}^{(1)},{\bf k}^{(1)}}
- Y'_{\e, {\bf h}^{(2)},{\bf k}^{(2)}} \|_{L^\infty} \\[3mm]
\qquad \qquad \qquad + r_\e^{2m-4} \,  \sup_{j=1, \ldots, n} \, \| (
\phi_{\e, {\bf h}^{(1)},{\bf k}^{(1)}} -
\phi_{\e, {\bf h}^{(2)},{\bf k}^{(2)}} )\,
|_{\bar B_{j ,2 r_\e} - B_{j, r_\e}} (r_\e \, \cdot)
 \|_{{\mathcal C}^{4, \alpha} (\bar B_{2} - B_{1} )} \qquad \qquad  \\[3mm]
\qquad \qquad \qquad  \leq  c_\kappa \, (r_\e^{2m - 3}  + \e^{2m-2}
\, r_\e^{2-2m-\delta}) \, \| ({\bf h}^{(1)}-{\bf h}^{(2)}, {\bf
k}^{(1)}-{\bf k}^{(2)})\|_{( {\mathcal C}^{4, \alpha})^{n} \times
({\mathcal C}^{2, \alpha})^{n} },
\end{array}
\label{eq:6.14}
\end{equation}
provided the components of ${\bf h}, {\bf h}^{(1)}, {\bf h}^{(2)}$,
${\bf k}, {\bf k}^{(1)}, {\bf k}^{(2)}$ satisfy (\ref{eq:NPR2}) and
(\ref{eq:NPR3}). \label{pr:6.22}
\end{prop}

The remaining of this section is devoted to the proof of this
result.

\begin{proof} To begin with, using the analysis of \S 6, we can expand~:
\[
{\bf s} (\omega + i \del \bar \del  f  ) =  {\bf s} (\omega)
-\tfrac{1}{2} \, P^*_\omega \, P_\omega \, f - \tfrac{1}{2} \, J
\, X_{\bf s} f + \tfrac{i}{2} \, X_{\bf s} \, f + Q_\omega
(\nabla^{2} \, f)
\]
The structure of the nonlinear operator $Q_\omega$ is quite
complicated but in each $\bar B_{j,\bar r_0}$, this operator enjoys
the following decomposition
\begin{equation}
\begin{array}{rllllll}
Q_\omega ( \nabla^{2} f ) & = &  \sum_{q} B_{q,4,2}(\nabla^{4}
f, \nabla^{2} f) \, C_{q,4,2} ( \nabla^{2} f) \\[3mm]
& + & \sum_{q} B_{q,3,3}(\nabla^{3} f, \nabla^{3} f) \,
C_{q,3,3} ( \nabla^{2} f ) \\[3mm]
& + & |z| \, \sum_{q} B_{q,3,2}(\nabla^{3} f, \nabla^{2}
\varphi) \, C_{q,3,2} ( \nabla^{2} f) \\[3mm]
& + & \sum_{q} B_{q,2,2}(\nabla^{2} f , \nabla^{2} f) \, C_{q,2,2} (
\nabla^{2} f)
\end{array}
\label{eq:6.3}
\end{equation}
where the sum over $q$ is finite, the operators $(U,V) \longmapsto
B_{q,a,b} (U,V)$ are bilinear in the entries and have coefficients
which are smooth functions on $\bar B_{j,\bar r_0}$. The nonlinear
operators $W \longmapsto C_{q,a,b} (W)$ have Taylor expansions (with
respect to $W$) whose coefficients are smooth functions on $\bar
B_{j,\bar r_0}$. These facts follow at once from the expression of
the scalar curvature of ${\bf s}( \omega_0 + i \, \del \, \bar \del
\, f)$ in local coordinates \cite{ap}.

The equation we would like to solve in $\bar M_{r_\e}$ reads
\[
{\bf s} (\omega  + i\del \bar \del  f ) -  \zeta - \mu_{\bf s} - \mu
= 0
\]
where
\[
\zeta = \ip{\xi_\omega, X_{\bf s} + X'} - \tfrac{1}{2} \, J \,
(X_{\bf s} + X') f - \tfrac{i}{2} \, (X_{\bf s} + X') \, f
\]
Observe that, using the analysis of section 5, we find that
\[
-d \zeta = (\omega + i \del \, \bar \del f)  \, (X_{\bf s} + X', -)
\]
Using the above expansion together with
(\ref{eq:solutioninitial}), we can rewrite this equation as
\[
- \tfrac{1}{2} \, P^*_\omega \, P_\omega \, f + i \, X_{\bf s} \, f
+ Q_\omega (\nabla^2 f) - \ip{\xi_\omega , X'} + \tfrac{i}{2} \, X'
\, f + \tfrac{1}{2} \, J\, X' \, f -\mu = 0
\]
Assuming that $X' \in \frh'$ and keeping in mind that we work
equivariantly so that the function $f$ is now assumed to be $K$
invariant, this simplifies into
\[
\tfrac{1}{2} \, P^*_\omega \, P_\omega  \, f + \ip{\xi_\omega , X'
} + \mu = \tfrac{1}{2} \, J \, X' \, f  +  Q_\omega (\nabla^{2} \,
f)
\]

We set
\[
\begin{array}{rlllllll}
f & : = & \e^{2m-2} \, \Gamma + W + \phi \\[3mm]
X' & : = & \e^{2m-2} \, Y' + Z' \\[3mm]
\mu & : = & \e^{2m-2} \, \lambda + \nu
\end{array}
\]
where $\Gamma$, $Y' \in \frh'$ and $\nu$ are defined in
(\ref{eq:gamma}), (\ref{eq:lambda}) and (\ref{eq:Xprime}) and $W =
W_{\e, {\bf h}, {\bf k}}$ is defined in (\ref{eq:doublev}).

It is easy to se that  the system of equations which remains to be
solved on $\bar M_{r_{\e}}$ can be formally written as
\[
\tfrac{1}{2}\, P^*_\omega \, P_\omega \phi + \ip{\xi_\omega , Z'} +
\nu  = {\mathcal Q} ( \phi , Z')
\] where the operator
${\mathcal Q} =  {\mathcal Q} (\e, {\bf h}, {\bf k} ; \cdot , \cdot
)$ is defined by
\[
\begin{array}{rllllll}
{\mathcal Q} (\e, {\bf h}, {\bf k} ; \phi, Z') & : = & \displaystyle
- \tfrac{1}{2} \, P^*_\omega \, P_\omega \, W + \tfrac{1}{2} \,
J\, (\e^{2m-2} \, Y' + Z' ) \, (\e^{2m-2} \, \Gamma  + W  +
\phi) \\[3mm]
&  + & Q_\omega (\nabla^{2} \, (\e^{2m-2} \, \Gamma  + W  + \phi ))
\end{array}
\]

We need the~:
\begin{defin}
Given $ \bar r \in (0, r_0)$, $\ell \in {\mathbb N}$, $\alpha \in
(0,1)$ and $\delta \in {\mathbb R}$, the weighted space ${\mathcal
C}^{\ell, \alpha}_{\delta} (\bar M_{\bar r })$ is defined to be the
space of functions $f \in {{\mathcal C}}^{\ell , \alpha} (\bar
M_{\bar r})$ endowed with the norm
\[
\| f \|_{{{\mathcal C}}^{\ell, \alpha}_{\delta} (\bar M_{\bar r})} :
= \| f \|_{{{\mathcal C}}^{\ell , \alpha} (\bar M_{r_0})} +
\sum_{j=1}^n \, \sup_{\bar r \leq r \leq r_0} r^{-\delta} \, \| f
|_{(\bar B_{j, r_{0}} -B_{j,\bar r})} (r \, \cdot) \|_{{{\mathcal C}
}^{\ell , \alpha} ( \bar B_{2} - B_{1})}
\]\label{de:6.1}
\end{defin}

Next, we consider an extension (linear) operator
\[
{\mathcal E}_{\bar r , \delta'} : {\mathcal C}^{0, \alpha}_{\delta'}
(\bar M_{\bar r}) \longrightarrow {\mathcal C}^{0, \alpha}_{\delta'}
(M^*) ,
\]
which is defined as follows~:
\begin{itemize}
\item[(i)]  In $M_{\bar r}$, ${\mathcal E}_{\bar r , \delta'} \, (f) = f$. \\

\item[(ii)] In each $\bar B_{j ,\bar r } - B_{j , \bar r/2}$
\[
{\mathcal E}_{\bar r , \delta'} \, (f) (z)  = \displaystyle
\tfrac{{2 \, |z| -\bar r  }}{{\bar r} }\,  f \left( \bar r \,
\tfrac{z}{{|z|}}
\right). \\
\]

\item[(iii)] In each $\bar B_{j , \bar  r/2}$,  $ {\mathcal E}_{\bar r , \delta'} \, (
f ) =  0 $.
\end{itemize}

It is easy to check that there exists a constant $c = c (\delta')
>0$, independent of ${\bar r} \in (0, r_0)$, such that
\begin{equation}
\| {\mathcal E}_{\bar r , \delta'} ( f ) \|_{{\mathcal C}^{0,
\alpha}_{\delta'} (M^{*})} \leq \, c \, \| f \|_{{\mathcal C}^{0,
\alpha}_{\delta'} (\bar M_{\bar r})} . \label{eq:6.8}
\end{equation}

The equation we would like to solve can now be written as
\begin{equation}
\tfrac{1}{2}\, P^*_\omega \, P_\omega  \, \phi + \ip{\xi_\omega ,
Y' } + \nu  + {\mathcal E}_{r_\e, \delta-4} \circ {\mathcal Q}_\e (
\phi , Y') =0 \label{pac-0701} \end{equation} We fix $\delta \in
(4-2m, 5-2m)$. If ${\mathcal G}_\delta$ denotes a right inverse for
${\mathcal L}_\delta$ which is provided by
Proposition~\ref{pr:f-5.1}, we just need to solve
\[
(\phi, Z' , \nu ) =  {\mathcal N} (\e, {\bf h}, {\bf k}, \phi , Z' )
\]
where  $\phi \in {\mathcal C}^{4, \alpha}_\delta (M^*)^K \oplus
{\mathcal D}$, $Z' \in {\mathfrak h}'$, $\nu \in {\mathbb R}$ and
where the nonlinear operator ${\mathcal N} (\e, {\bf h}, {\bf k} ;
\, \cdot \, , \, \cdot \, )$ is defined by
\[
 {\mathcal N} (\e, {\bf h}, {\bf k} ; \, \cdot \, , \,
\cdot \, ) : = {\mathcal G}_\delta \circ {\mathcal E}_{r_\e,
\delta-4} \circ {\mathcal Q}
\]
We set
\[
{\mathcal F} : = ( {\mathcal C}^{4, \alpha}_{\delta} (\bar M^*)^K
\oplus {\mathcal D} ) \times {\frh'} \times {\mathbb R}
\]
which is endowed with the product norm.

The existence of a solution to this system depends is a simple
consequence of the fixed point theorem for contraction mappings once
the following estimates are proved.
\begin{lemma}
There exists $c  >0$, $c_\kappa  > 0$ and there exists $\e_\kappa =
\e (\kappa) >0$ such that, for all $\e \in (0, \e_\kappa)$
\begin{equation}
\| {\mathcal N}  (\e , {\bf h},{\bf k} ; 0 , 0) \|_{{\mathcal F}}
\leq c_\kappa \, (r_\e^{2m+1} + \e^{4m-4} \, r_\e^{6 -4m -\delta }),
\label{eq:6.1000}
\end{equation}
and
\begin{equation}
\begin{array}{llllll}
\| {\mathcal N}  (\e , {\bf h},{\bf k} ; \phi^{(1)} , Z' \, ^{(1)} )
- {\mathcal N} (\e , {\bf h},{\bf k} ; \phi^{(2)} , Z' \, ^{(2)})
\|_{{\mathcal F}} \\[3mm]
\hspace{40mm}\leq c_\kappa \, \e^{2m-2} \, r_\e^{6-4m-\delta} \, \|
(\phi^{(1)} - \phi^{(2)} , Z' \, ^{(1)} - Z' \, ^{(2)}, 0)
\|_{{\mathcal F}}
\end{array}\label{eq:6.1001}
\end{equation}
Finally,
\begin{equation}
\begin{array}{llllll}
\| {\mathcal N}  (\e , {\bf h}^{(1)},{\bf k}^{(1)} ; \phi, Z' ) -
{\mathcal N} (\e , {\bf a} , {\bf h}^{(2)},{\bf k}^{(2)} ; \phi , Z'
) \|_{{\mathcal F}} \\[3mm]
\hspace{30mm}\leq c_\kappa \, (r_\e^{2m-3} + \e^{2m-2} \,
r_\e^{2-2m-\delta}) \, \| ({\bf h}^{(1)} - {\bf h}^{(2)}, {\bf
k}^{(1)} - {\bf k}^{(2)}) \|_{( {\mathcal C}^{4, \alpha})^n \times
({\mathcal C}^{2, \alpha} )^n}
\end{array}
\label{eq:6.1002}
\end{equation}
provided $(\phi, Z',0), (\phi^{(1)}, Z' \, ^{(1)},0) , (\phi^{(2)},
Z' \, ^{(2)},0) \in {\mathcal F}$ satisfy
\[
\| (\phi  , Z' ,0) \|_{{\mathcal F}} + \| (\phi^{(1)} , Z' \, ^{(1)}
,0) \|_{{\mathcal F}} +\| (\phi^{(2)} , Z' \, ^{(2)} ,0)
\|_{{\mathcal F}} \leq \, 6 \, c_\kappa \, (r_\e^{2m+1} +  \e^{4m-4}
\, r_\e^{6 -4m -\delta }) ,
\]
and the components of ${\bf h}, {\bf h}^{(1)}, {\bf h}^{(2)}$, ${\bf
k}, {\bf k}^{(1)}, {\bf k}^{(2)}$ satisfy (\ref{eq:NPR2}) and
(\ref{eq:NPR3}). \label{le:f-8.1}
\end{lemma}
\begin{proof} The proof of these estimates can be follows what is already
done in \cite{ap} and \cite{ap2} with minor modifications. We
briefly recall how the proof of the first estimate is obtained and
leave the proof of the second and third estimates to the reader.
First, we use the result of Proposition~\ref{pr:f-5.5} to estimate
\begin{equation}
\| W \|_{{\mathcal C}^{4, \alpha}_{3-2m} (\bar M_{r_\e})} \leq
c_\kappa \, r_\e^{2m+1} . \label{eq:jd}
\end{equation}
Now observe that, by construction, $\Delta^{2} \, W  = 0$ in each
$\bar B_{r_0/2} - B_{r_\e}$ (here $\Delta$ is the Euclidean
Laplacian), hence
\[
L_\omega \, W  = (L_\omega + \tfrac{1}{2} \, \Delta^{2} ) \, W
\]
in this set. Making use of the fact that the coordinates near $p_j$
are chosen to be normal, we get the existence of a constant
$c_\kappa > 0$ such that
\[
\| L_\omega \, W  \|_{{\mathcal C}^{0, \alpha}_{\delta -4} ( \bar
M_{r_\e})} \leq c_\kappa \, r_\e^{2m+1}.
\]
Note that this is where we implicitly use the fact that $\delta <
5-2m$.

Next, we estimate the norm of $\e^{2m-2} \, J \, Y'  \, ( \e^{2m-2}
\, \Gamma  + W )$ by
\[
\| \e^{2m-2} \, J \, Y' \, ( \e^{2m-2} \, \Gamma  + W )
\|_{{\mathcal C}^{0, \alpha}_{\delta -4} ( \bar M_{r_\e})} \leq c \,
\e^{4m -4}.
\]
where the constant $c_\delta$ does not depend on $\kappa$ provided
$\e$ is chosen small enough, say $\e \in (0, \e_\kappa)$.

Finally, we use the structure of the nonlinear operator $Q_\omega$
as described in (\ref{eq:6.3}) together with the estimate
(\ref{eq:jd}) to get
\[
\| Q_\omega ( \nabla^{2}  ( \e^{2m-2} \, \Gamma + W)) \|_{{\mathcal
C}^{0, \alpha}_{\delta-4} (\bar M_{r_\e})} \leq c \, \e^{4m-4} \,
r_\e^{6-4m-\delta}
\]
for some constant $c_\delta >0$ which does not depend on $\kappa$
provided $\e$ stays small enough. The first estimate then follows at
once. \end{proof}

Reducing $\e_\kappa >0$ if necessary, we can assume that,
\begin{equation}
c_\kappa \, \e^{2m-2} \, r_\e^{6-4m-\delta}  \leq \tfrac{1}{2}
\label{eq:6.1400}
\end{equation}
for all $\e \in (0, \e_\kappa )$. Then, the estimates
(\ref{eq:6.1000}) and (\ref{eq:6.1001}) in the above Lemma are
enough to show that
\[
(\phi, X' , \nu ) \longmapsto {\mathcal N} (\e, {\bf h} , {\bf k} ;
\phi, X' )
\]
is a contraction from
\[
\{ (\phi, X' , \nu ) \in {\mathcal F}  \quad : \newline \quad
\|(\phi , X' , \nu )\|_{{\mathcal F}} \leq 2 \,  c_\kappa \, (
r_\e^{2m+1} + \e^{4m-4} \, r_\e^{6 -4 m -\delta } ) \} ,
\]
into itself and hence has a unique fixed point $(\phi_{\e, {\bf h},
{\bf k}}, Y'_{\e, {\bf h}, {\bf k}} , \nu_{\e, {\bf h}, {\bf k}})$
in this set. This fixed point yields a solution of (\ref{pac-0701})
in $\bar M_{r_\e}$ and hence provides an extremal \K\ form on $\bar
M_{r_\e}$. The estimates in Proposition~\ref{pr:6.22} follow at once
from the estimates in Lemma~\ref{le:f-8.1}, increasing the value of
$c_\kappa$ and reducing $\e_\kappa$ if this is necessary.
\end{proof}

\subsection{Perturbation of $\eta$}

We perform an analysis similar to the one we have done in the
previous section starting with $\tilde {\mathbb C}^m$, the blow up
of ${\mathbb C}^m$ at $0$, endowed with the Burns-Simanca's metric $\eta$.

Given $ a >0$, we now consider on $\tilde {\mathbb C}^m$, the
perturbed \K\ form
\[
\tilde \eta = a^2 \, \eta + i \, \del \, \bar \del \, f ,
\]
Everything we will do will be uniform in $a$ as long as this
parameter remains both bounded from above and bounded away from $0$.
In fact, we will apply our analysis when $a$ satisfies
\begin{equation}
a_{min}  \leq a \leq   a_{max}, \label{eq:aaa}
\end{equation}

Given $\nu \in {\mathbb R}$ and a $K$-invariant killing field $X \in
\frk$, we would like to solve the equation
\begin{equation} {\bf s} \, (a^{2} \, \eta + i \, \del \, \bar \del
\, f) = \e^{2} \, \zeta \label{eq:6.17}
\end{equation} in $N_{R_\e/a}$ where
\[
\zeta = \tilde \zeta  - \tfrac{1}{2} \, J \, X \, f - \tfrac{i}{2}
\, X \, f
\]
and $\tilde \zeta$ is the solution of
\[
- d\tilde \zeta =  a^{2} \, \eta ( X , - )
\]
normalized so that the mean of $\zeta$ over $\del B_{R_\e/a}$ is
prescribed by
\[
\tfrac{1}{{|\del B_{R_\e/a}|}} \, \int_{\del B_{R_\e/a}} \zeta =
\nu
\]
Observe that, using the analysis of section 5, we find that
\[
- d \zeta = (a^2 \, \eta + i \, \del \, \bar \del f)  \, (X, -)
\]
Also, working with $K$-invariant functions, since $X \in \frk$, we
have $X\, f =0$.

It will be convenient to denote
\[
\bar N_R :  =  \tilde {\mathbb C}^m  - C_R
\]
We fix
\begin{equation}
R_\e : =  \tfrac{r_\e}{\e} \label{eq:RRe}
\end{equation}
Assume we are given $h \in {\mathcal C}^{4, \alpha}(\del B_1)^K$ and
$k \in {\mathcal C}^{2, \alpha}(\del B_1)^K$ satisfying
\begin{equation}
\| h \|_{{\mathcal C}^{4, \alpha} (\del B_1)} + \| k \|_{{\mathcal
C}^{2, \alpha} (\del B_1)} \leq \kappa \, R_\e^{3-2m},
\label{eq:f-23}
\end{equation}
satisfying
\[
\int_{\del B^1} (4mh-k) =0
\]
where $\kappa >0$ will be fixed later on. We define in $N_{R_\e/a}$
the function
\begin{equation}
W_{\e, a , h,k} : = \tilde \chi \, (W^i_{h , k } (a \, \cdot / R_\e)
, \label{eq:f-24}
\end{equation}
where $\tilde \chi$ is a cutoff function which is identically equal
to $1$ in $C_2$ and identically equal to $0$ in $N_1$ and where
$W^{i}_{h,k}$ has been defined in Proposition~\ref{pr:f-5.5}.

We will assume that the parameter $\nu$ and the vector field $X$ are
uniformly bounded in $N_{R_\e/a}$. To be more specific, we agree
that
\[
|\nu| + \e^{-1} \, \| X \|_{L^\infty (N_{2R_\e/a})} \leq c
\]

Following the analysis already performed in the previous section, we
have the~:
\begin{prop} Given $\delta \in (0,1)$, there exist $c
>0$ (independent of $\kappa$) and $\e_\kappa = \e (\kappa) > 0$ such
that, for all $\e \in (0, \e_\kappa)$  there exists a function
$\phi_{\e, a, \nu, X, h,k} \in {\mathcal C}^{4, \alpha} (\bar
N_{R_\e/a})$ such that the scalar curvature of the \K\ form
\[
\eta_{\e, a ,\nu, X , h , k  } : =  a^{2} \, \eta  + i \, \del \,
\bar \del \, (W_{\e, a , h , k } + \phi_{\e, a , \nu , X, h , k } ),
\]
defined on $\bar N_{R_\e/a}$, is given by
\[
- d {\bf s} (\eta_{\e, a ,\nu, X , h , k  } ) = \e^2 \, \eta_{\e, a
,\nu, X , h , k  } (X, -)
\]
with
\[
\tfrac{1}{{|\del B_{R_\e/a}|}} \, \int_{\del B_{R_\e/a}} {\bf s}
(\eta_{\e, a ,\nu, X , h , k  } )   = \e^2
 \, \nu
\] Moreover
\[
\| \phi_{\e, a , \nu , X,  h , k } (R_\e \, \cdot /a ) \|_{{\mathcal
C}^{4, \alpha} (\bar B_{1} -B_{1/2})} \leq c \, R_\e^{3-2m} ,
\]
for some constant $c >0$ independent of $\kappa$. In addition, we
have
\begin{equation}
\begin{array}{llll} \| \phi_{\e, a , \nu , X,  h , k } (R_\e \,
\cdot /a )  - \phi_{\e, a , \nu' , X',  h' , k' } (R_\e \, \cdot /a'
)  \|_{{\mathcal C}^{4, \alpha} ( \bar B_{1} - B_{1/2} )}  \\[3mm]
\hspace{20mm}  \leq c_\kappa \, ( R_\e^{\delta -1} \, \| (h - h', k
- k' )\|_{{\mathcal C}^{4, \alpha} \times {\mathcal C}^{2, \alpha} }
+ R_\e^{3-2m} \,( | \nu - \nu'|  +  \| X-X' \|_{L^\infty} + |a-a'|)
).
\end{array}
\label{eq:6.29}
\end{equation}
\label{pr:6.2}
\end{prop}

Again the rest of the section is devoted to the proof of this
result.

\begin{proof}
Using the fact that
\[
{\bf s}  ( a^{2} \, \eta + i \, \del \, \bar \del \, f ) = {\bf s} (
a^{2} \, (\eta + i \, a^{-2} \, \del \, \bar \del \, f  )) = a^{- 2}
\, {\bf s}  ( \eta + i \, a^{- 2} \, \del \, \bar \del \, f )
\]
We see that the scalar curvature of $\tilde \eta$ can be expanded as
\begin{equation}
{\bf s} (a^{2} \, \eta + i \, \del \, \bar \del \, f) =  -
\tfrac{1}{2} \, a^{-2} \, P_{\eta}^* \, P_\eta \, f + a^{- 2} \,
Q_{\eta} (a^{- 2} \, \nabla^{2} f ) , \label{eq:6.1666}
\end{equation}
since the scalar curvature of $\eta$ is identically equal to $0$.
Again, the structure of the nonlinear operator $Q_{\eta}$ is quite
complicated but away from the exceptional divisor, it enjoys a
decomposition similar to the one described in the previous section.
Indeed, we know from \cite{ap2} that we can decompose
\[
\begin{array}{rlllll}
Q_{\eta} ( \nabla^{2} f) & = & \sum_{q} B_{q,4,2}(\nabla^{4}
f, \nabla^{2} f) \, C_{q,4,2} ( \nabla^{2} f) \\[3mm]
& + & \sum_{q} B_{q,3,3}(\nabla^{3} f, \nabla^{3} f) \,
C_{q,3,3} ( \nabla^{2} f) \\[3mm]
& + & \sum_{q} |u|^{1-2m} \, B_{q,3,2}(\nabla^{3} f, \nabla^{2} f)
\, C_{q,3,2} ( \nabla^{2}f)
\\[3mm]
& +& \sum_{q} |u|^{-2m} \, B_{q,2,2}(\nabla^{2} f, \nabla^{2} f) \,
C_{q,2,2} ( \nabla^{2} f)
\end{array}
\]
where the sum over $q$ is finite, the operators $(U,V)
\longrightarrow B_{q,j,j'} (U,V)$ are bilinear in the entries and
have coefficients which are bounded functions in ${\mathcal C}^{0,
\alpha} (\bar C_1)$. The nonlinear operators $W \longrightarrow
C_{q,a,b} (W)$ have Taylor expansion (with respect to $W$) whose
coefficients are bounded functions on ${\mathcal C}^{0, \alpha}
(\bar C_1)$.

We set
\[
L_\eta : = - \tfrac{1}{2} \, P^*_\eta \, P_\eta
\]
Replacing in (\ref{eq:6.1666}) the function
\[
f = W  + \phi,
\]
where $W = W_{\e, h,k}$ has been defined in (\ref{eq:f-24}), we see
that (\ref{eq:6.17}) can be written as
\begin{equation}
L_{\eta} \, \phi  + Q_{\eta} ( a^{-2} \, \nabla^{2} ( W + \phi ) ) =
\e^{2} \, a^{2} \, \zeta , \label{eq:6.99}
\end{equation}
which we would like to solve in $\bar N_{R_\e/a}$. Remember that the
function $\zeta$ solves
\[
- d \zeta =  ( a^{2} \, \eta + i \, \del \, \bar \del \, ( W+ \phi))
( X , - )
\]
with
\[
\tfrac{1}{{|\del B_{R_\e/a}|}} \, \int_{\del B_{R_\e/a}} \zeta =
\nu
\]

We will need the~:
\begin{defin}
Given $ \bar R > 1 $, $\ell \in {\mathbb N}$, $\alpha \in (0,1)$ and
$\delta \in {\mathbb R}$, the weighted space ${\mathcal C}^{\ell,
\alpha}_{\delta} (\bar N_{\bar R})$ is defined to be the space of
functions $f  \in {{\mathcal C}}^{\ell, \alpha} (\bar N_{\bar R})$
endowed with the norm
\[
\| f \|_{{{\mathcal C}}^{\ell, \alpha}_{\delta} (\bar N_{\bar R})} :
= \| f \|_{{{\mathcal C}}^{\ell, \alpha} (\bar N_1)} + \sup_{1 \leq
R \leq \bar R}  R^{-\delta} \, \| f (R \, \cdot) \|_{{{\mathcal C}
}^{\ell , \alpha} ( \bar B_{1} - B_{1/2})}
\]\label{de:6.2}
\end{defin}

For each $\bar R \geq 1$, will be convenient to define an
"extension" (linear) operator
\[
\tilde {\mathcal E}_{\bar R} : {\mathcal C}^{0, \alpha}_{\delta'}
(N_{\bar R}) \longrightarrow {\mathcal C}^{0, \alpha}_{\delta'}
(\tilde {\mathbb C}^m ) ,
\]
as follows~:
\begin{itemize}

\item[(i)] In $\bar N_{\bar R}$,  $\tilde {\mathcal E}_{\bar R} \, ( f ) =
f$,\\

\item[(ii)] in $C_{2\bar R} - C_{\bar R} $
\[
\tilde {\mathcal E}_{\bar R} \, (f ) (u) = \tfrac{{2 \, {\bar R} -
|u|}}{{\bar R}} \, f \left(\bar R \, \tfrac{u}{{|u|}}\right),
\]

\item[(iii)] in $C_{2 \, \bar R }$, $\tilde {\mathcal E}_{\bar R} \,
( f )=0$.

\end{itemize}
It is easy to check that there exists a constant $c = c( \delta' )
>0$, independent of $\bar R \geq 2$, such that
\begin{equation}
\| \tilde {\mathcal E}_{\bar R} ( f ) \|_{{\mathcal C}^{0,
\alpha}_{\delta'} (\tilde {\mathbb C}^m)} \leq \, c \, \| f
\|_{{\mathcal C}^{0, \alpha}_{\delta'} (N_{\bar R})} ,
\label{eq:6.20}
\end{equation}
and furthermore, one can arrange easily for $\tilde {\mathcal
E}_{\bar R}$ to depend smoothly on $\bar R$.

The equation we will solve can be rewritten as
\begin{equation}
L_{\eta} \, \phi   = \tilde {\mathcal E}_{R_\e/a} \left( \left(
\e^{2} \, a^{2} \, \zeta  - Q_{\eta} (a^{-2}  \, \nabla^{2} ( W +
\phi )) - L_{\eta} \, W \right) \right). \label{eq:f-25}
\end{equation}
We fix $\delta \in (0, 1)$ and use the result of
Proposition~\ref{pr:MP3}. This provides a right inverse $\tilde
{\mathcal G}_{\delta}$ for the operator $L_{\eta}$.

We can now rephrase the solvability of (\ref{eq:f-25}) as a fixed
point problem.
\begin{equation}
\phi = \tilde {\mathcal N} (\e, a, \nu, X, h, k  ; \phi )
\label{eq:6.22}
\end{equation}
where the nonlinear operator ${\tilde {\mathcal N}}$ is defined by
\[
\tilde {\mathcal N} (\e , a, \nu, X , h, k ; \phi ) : = \tilde
G_{\delta} \, \circ \tilde {\mathcal E}_{R_\e/a} \left( ( Q_{\eta} (
a^{-2} \, \nabla^{2} (W  + \phi )) ) - {\mathbb L}_{g_1} \, \tilde
H_{h, k} - \e^{2} \, a^2 \, (\nu +\zeta) \right)
\]
To keep notations short, it will be convenient to define
\[
\tilde {\mathcal F} : = {\mathcal C}^{4, \alpha}_\delta (\tilde
{\mathbb C}^m)
\]

The existence of a fixed point for ${\tilde {\mathcal N}}$ will
follow from the~:
\begin{lemma}
There exists $c >0$ (independent of $\kappa$),  $c_\kappa  >0$ and
there exists $\e_\kappa = \e (\kappa) >0$ such that, for all $\e \in
(0, \e_\kappa)$
\begin{equation}
\| \tilde {\mathcal N} ( \e, a, \nu, X,  h, k ; 0 ) \|_{\tilde
{\mathcal F}} \leq c \, R_\e^{3-2m -\delta} , \label{eq:6.23}
\end{equation} Moreover, we have
\begin{equation}
\| \tilde {\mathcal N} ( \e, a, \nu, X , h, k ; \phi ) - \tilde
{\mathcal N} ( \e, a , \nu, X, h, k  ; \phi' ) \|_{\tilde {\mathcal
F}} \leq c_\kappa \, R_\e^{3-2m -\delta} \, \| \phi - \phi'
\|_{\tilde {\mathcal F}} \label{eq:6.24}
\end{equation} and
\begin{equation}
\begin{array}{rllllll}
\| \tilde {\mathcal N} ( \e, a , \nu, X,  h, k  ; \phi ) - \tilde
{\mathcal N} ( \e, a' , \nu' , X',  h', k' ; \phi ) \|_{\tilde
{\mathcal F}} \leq   c_\kappa \, ( R_\e^{-1} \, \| (h-h',
k-k')\|_{\mathcal C^{4, \alpha} \times {\mathcal C}^{2, \alpha}} \\[3mm]
\hspace{40mm} +   R_\e^{3-2m -\delta} \, ( |\nu' - \nu| + \| X
-X'\|_{L^\infty}+ |a' - a| ) )
\end{array}
\label{eq:6.25}
\end{equation} provided $\phi , \phi' \in \tilde {\mathcal F}$, satisfy
\[
\|\phi \|_{\tilde {\mathcal F}} + \|\phi'\|_{\tilde {\mathcal F}}
\leq 4 \, c \, R_\e^{3-2m -\delta} ,
\]
and $h, h'$ and $k,k'$ satisfy (\ref{eq:f-23}). \label{le:f-8.3}
\end{lemma}
The proof of these estimates being identical to the one in
\cite{ap2}, we omit it.

Reducing $\e_\kappa >0$ if necessary, we can assume that,
\begin{equation}
c_\kappa \, R_\e^{3-2m-\delta}  \leq \tfrac{1}{2} \label{eq:6.28}
\end{equation} for all $\e \in (0, \e_\kappa )$. Then, the estimates
(\ref{eq:6.23}) and (\ref{eq:6.24}) in the above Lemma are enough to
show that
\[
\phi \longmapsto \tilde{\mathcal N} (\e, a , \nu, X, h , k ; \phi )
\]
is a contraction from
\[ \{ \phi  \in \tilde {\mathcal F} \quad : \quad \|
\phi \|_{\tilde {\mathcal F}} \leq 2 \, c \, R_\e^{3-2m-\delta} \} ,
\]
into itself and hence has a unique fixed point $\phi_{ \e, a , \nu,
X, h, k}$ in this set. This fixed point is a solution of
(\ref{eq:6.99}) in $\bar N_{R_\e}$ and hence provides a constant
scalar curvature \K\ form on $\bar N_{R_\e}$. The estimates in
Proposition~\ref{pr:6.2} follow at once from the estimates in
Lemma~\ref{le:f-8.3}, increasing the value of $c_\kappa$ and
reducing $\e_\kappa$ if this is necessary. \end{proof}

\section{Gluing the pieces together}

Building on the analysis of the previous sections we complete the
proof of Theorem~\ref{mainthm}. As far as technicalities are
concerned the proof is identical to the one in \cite{ap2} therefore,
we shall only emphasize the differences in the present framework.

Before we proceed, a word about notations. In this section
${\mathcal O}_{{\mathcal C}^{\ell, \alpha}} (A)$ will refer to a
function whose ${\mathcal C}^{\ell, \alpha}$-norm is bounded by $A$
times a constant independent of $\e$ and also independent of
$\kappa$ provided $\e$ is chosen small enough (but which might
depend on $m$, $\omega$, the points $p_j$ and the coefficients
$a_j$). In general this function will be a nonlinear operator of the
data.

We first exploit the result of Proposition~\ref{pr:6.22}. Given
boundary data
\[
{\bf h} : = (h_1, \ldots, h_n), \qquad {\bf k} : = (k_1, \ldots,
k_n)
\]
so that
\[
\int_{\del B_1} k_j =0
\]
we can apply the result o Proposition~\ref{pr:6.22} to define on
$\bar M_{r_\e}$ a \K\ form $\omega_{\e, {\bf h}, {\bf k}}$  which
can be written as
\[
\omega_{\e, {\bf h}, {\bf k}} = i \, \del \, \bar \del
(\tfrac{1}{2} \, |z|^2 + \varphi^j(z) + \e^{2m-2} \, \Gamma_{\e,
{\bf h}, {\bf k}} + W_{\e, {\bf h}, {\bf k}} + \phi_{\e, {\bf h},
{\bf k}} )
\]
in $\bar B_{j, 2r_\e} -B_{j, r_\e}$ where the function $\varphi^j
(z) = {\mathcal O} (|z|^4)$ is the one which appears in
Proposition~\ref{masterpiece} so that
\[
\omega = i \, \del \, \bar \del (\tfrac{1}{2} \, |z|^2 +
\varphi^j(z))
\]
near $p_j$. We define the function
\[
\psi^{o,j} : = \left( \varphi^j + \e^{2m-2} \, \Gamma_{\e, {\bf h},
{\bf k}} + W_{\e, {\bf h}, {\bf k}} + \phi_{\e, {\bf h}, {\bf
k}}\right) (r_\e \, \cdot )
\]
in $\bar B_2 - B_1$. Collecting the result of
Proposition~\ref{pr:6.22}, the definition of $W^o_{\e, {\bf h}, {\bf
k}}$ given in (\ref{eq:doublev}) and the expansion of $\Gamma$ given
in Lemma~\ref{le:NPR1} we find that the function $\psi^{o,j}$ can be
expanded as
\[
\psi^{o,j}  = -  \tfrac{1}{{m-2}} \, a_j \, \e^{2m-2} \,
r_\e^{4-2m} \, |\cdot|^{4-2m} + W^o_{h_j, k_j} + {\mathcal
O}_{{\mathcal C}^{4, \alpha}} (r_\e^{4})
\]
in dimension $m\geq 3$ while, in dimension $m=2$, in view of the
expansion of $\Gamma$ Lemma~\ref{le:NPR1}, we have
\[
\psi^{o,j} - \e^{2} \, (a_j \, \log r_\e + \e^2 \, b_j)  = a_j \,
\e^{2} \, \log |\cdot| + W^o_{h_j, k_j} + {\mathcal O}_{{\mathcal
C}^{4, \alpha}} (r_\e^{4}).
\]
We will replace $\psi^{o,j}$ by $\psi_{o,j} - \e^2  \, ( a_j \, \log
r_\e + b_j)$ and there is no loss of generality in doing so since
changing the potential by some constant function does not alter the
corresponding \K\ forms.

According to Proposition~\ref{pr:6.22}, the scalar curvature of the
\K\ form $\omega_{\e, {\bf h}, {\bf k}}$ is given by
\[
{\bf s} (\omega_{\e, {\bf h}, {\bf k}})  = {\bf s}(\omega) +
\ip{\xi_{\omega_{\e, {\bf h}, {\bf k}}} , X_{\e, {\bf h}, {\bf k}}}
+ \mu_{\bf s} + \e^{2m-2} \, \lambda + \lambda_{\e, {\bf h}, {\bf
k}}
\]
where we have defined
\begin{equation}
X_{\e, {\bf h}, {\bf k}} : = X_{\bf s} + \e^{2m-2} \, Y' + Y'_{\e,
{\bf h}, {\bf k}} \in \frh' \label{hvf}
\end{equation}

We now exploit the result of Proposition~\ref{pr:6.2}. We choose
boundary data
\[
{\bf \tilde h} : = (\tilde h_1 , \ldots, \tilde h_n),  \qquad {\bf
\tilde k} : = (\tilde k_1 , \ldots, \tilde k_n)
\]
whose components satisfy
\begin{equation}
\int_{\del B_1} (4m \tilde h_j - \tilde k_j) = 0 \label{eq:concon}
\end{equation} as well as real positive parameters $\tilde {\bf a} : =
(\tilde a_1, \ldots, \tilde a_n)$. For each $j =1, \ldots, n$, we
apply the result of Proposition~\ref{pr:6.2} and define on $\bar
N_{R_\e /\hat a_j}$ the \K\ form
\[
\e^2 \, \eta_{\e, \hat a_j , \nu_j, X_j , \tilde h_j, \tilde k_j}
\]
where
\[
\hat a_j : = \tilde a_j^{\tfrac{1}{{2(m-1)}}}
\]
and
\[
\nu_j : = \tfrac{1}{{|\del B_{j, r\e}|}} \, \int_{\del B_{j, r_\e}}
{\bf s} (\omega_{\e, {\bf h},{\bf k}}) (r_\e \cdot)
\]
Moreover the scalar curvature of $\e^2 \, \eta_{\e, \hat a_j ,
\nu_j, X_j , \tilde h_j, \tilde k_j}$ satisfies
\[
- d \, {\bf s} (\e^2 \, \eta_{\e, \hat a_j , \nu_j, X_j , \tilde
h_j, \tilde k_j})  = \eta_{\e, \hat a_j , \nu_j, X_j , \tilde h_j,
\tilde k_j} ( X_j, -)
\]
with
\begin{equation}
\int_{\del B_1}  \zeta (R_\e \cdot / \hat a_j) = \nu_j
\label{eq:nuj}
\end{equation}
and $X_j$ is the lift to $\bar N_{R_e/\hat a_j}$ of the holomorphic
vector field $X_{\e, {\bf h}, {\bf k}}$ defined in $B_{j, r_\e}$. As
explained at the end of section 7, this lifting can be performed in
normal  $K$-linear coordinates so that the vector field $X_j$ is a
$K$-invariant Killing vector field for the metric $\hat a_j^{2} \,
\eta$.

According to the analysis of section 6 and the result of
Proposition~\ref{pr:6.2}, the \K\ form $\eta$ can be written as
\[
\e^2 \, \eta_{\e, \hat a_j , \nu_j, X_j , \tilde h_j, \tilde k_j} =
i \, \del \bar \del \left( \e^2 \, \hat a_j^{2} \, E_m (u) + \e^2 \,
W_{\e, \hat a_j, \tilde h_j, \tilde k_j} + \e^2 \, \phi_{\e, \hat
a_j, \nu_j, X_j, \tilde h_j, \tilde k_j} \right)
\]
in $B_{R_\e/\hat a_j} - B_{R_\e/2\hat a_j}$. We define the function
\[
\psi^{i,j}  : = \left( \e^2 \, \hat a_j^{2} \, (E_m  - \tfrac{1}{2}
\, |\cdot|^2 ) + \e^2 \, W_{\e, \hat a_j, \tilde h_j, \tilde k_j} +
\e^2 \, \phi_{\e, \hat a_j, \nu_j, X_j, \tilde h_j, \tilde k_j}
\right)  ( R_\e \cdot / \hat a_j),
\]
defined in $\bar B_{1} - B_{1/2}$. Using the analysis of section 6
as well as the result of Proposition~\ref{pr:6.2}, we have the
expansion
\[
\psi^{i,j}  = - \, \tilde a_j \, \e^{2m-2} \, r_\e^{4-2m} \, |\cdot
|^{4-2m} + \e^2 \, W^i_{\tilde h_j, \tilde k_j}  + {\mathcal
O}_{{\mathcal C}^{4, \alpha}}(r_\e^{4})
\]
in $\bar B_1 - B_{1/2}$, when $m \geq 3$ while we have
\[
\psi^{i,j} - \tilde a_j \, \e^2 \log \e  = \tilde a_j \, \e^{2} \,
\log |\cdot | + \e^2 \, W^i_{\tilde h_j, \tilde k_j}  + {\mathcal
O}_{{\mathcal C}^{4, \alpha}}(r_\e^{4})
\]
when $m=2$. Again, in dimension $m=2$ we will replace $\psi^{i,j}$
by $\psi^{i,j} - \tilde a_j \, \e^2 \log \e $ since this does not
affect the definition of the corresponding \K\ metric.

The proof now follows {\it verbatim} the proof in \cite{ap2}. We
first describe the connected sum construction. By construction,
\[
M_\e : = M \sqcup _{{p_{1}, \e}} N_1 \sqcup_{{p_{2},\e}} \dots
\sqcup _{{p_n, \e}} N_n ,
\]
is obtained by connecting  $M_{r_\e}$ with the truncated spaces
$N_{R_\e/ \tilde a_1}, \ldots, N_{R_\e/ \tilde a_n}$. The
identification of the boundary $\del B_{j , r_\e}$ in $M_{r_\e}$
with the boundary $\del N_{R_\e/\tilde a_j}$ of $N_{R_\e/ \tilde
a_j}$ is performed using the change of variables
\[
(z^{1} , \ldots, z^{m} )  = \e \, \hat a_j \, (u^{1} , \ldots,
u^{m}) ,
\]
where $(z^{1}, \ldots, z^{m} )$ are the coordinates in $B_{j , r_0}$
and $(u^{1}, \ldots, u^{m})$ are the coordinates in $C_1$.

The problem is now to determine these boundary data and parameters
in such a way that, for each $j=1, \ldots, n$ the functions
$\psi^{o,j} $ and $\psi^{i,j} $ have  have their partial derivatives
up to order $3$ which coincide on $\del B_{1}$.

In fact, we shall solve the following system of equations
\begin{equation}
\psi^{o,j} =\psi^{i,j} , \qquad \del_r \, \psi^{o,j} = \del_r \,
\psi^{i,j} , \qquad \Delta \, \psi^{o,j} = \Delta \, \psi^{i,j},
\qquad \del_r \, \Delta  \, \psi^{o,j} = \del_r \, \Delta \,
\psi^{i,j}, \label{eq:6.300}
\end{equation}
on $\del B_{1}$ where $r =|v|$ and $v= (v^{1}, \ldots, v^{m})$ are
coordinates in ${\mathbb C}^{m}$.

Let us assume that we have already solved this problem. The first
identity in (\ref{eq:6.300}) implies that $\psi^{o,j}$ and
$\psi^{i,j}$ as well as all their $k$-th order partial derivatives
with respect any vector field tangent to $\del B_{1}$, with $k \leq
4$, agree on $\del B_{1}$. The second identity in (\ref{eq:6.300})
then shows that $\del_r \psi^{o,j}$ and $\del_r \psi^{i,j}$ as well
as all their $k$-th order partial derivatives with respect any
vector field tangent to $\del B_{1}$, with $k \leq 3$, agree on
$\del B_{1}$. Using the decomposition of the Laplacian in polar
coordinates, it is easy to check that the third identity implies
that $\del_r^{2} \psi^{o,j}$ and $\del_r^{2} \psi^{i,j}$ as well as
all their $k$-th order partial derivatives with respect any vector
field tangent to $\del B_{1}$, with $k \leq 2$, agree on $\del
B_{1}$. And finally, the last identity in (\ref{eq:6.300}) implies
that $\del_r^{3} \psi^{o,j}$ and $\del_r^{3} \psi^{i,j}$ as well as
all their first order partial derivative with respect any vector
field tangent to $\del B_{1}$, agree on $\del B_{1}$.

Moreover, the scalar curvature of the \K\ form
\[
\omega^{o,j} : = i \, \del \, \bar \del \, ( \tfrac{1}{2} \,
|v|^{2} + \psi^{o,j} ),
\]
defined in $\bar B_{2} - B_{1}$ and the scalar curvature of the \K\
form
\[
\omega^{i,j} : = i \, \del \, \bar \del \, ( \tfrac{1}{2} \,
|v|^{2} + \psi^{i,j}),
\]
defined in $\bar B_{1} - B_{1/2}$, match on $\del B_1$ to produce a
${\mathcal C}^{2}$ function on $\bar B_2-B_{1/2}$. To see this
observe that both scalar curvature functions have the same mean
value on $\del B_1$ (this was precisely the purpose of
(\ref{eq:nuj})) and they satisfy
\[
-d {\bf s} =  \hat \omega (X , -)
\]
for the same vector field $X$. Since the \K\ form is already
${\mathcal C}^1$, we find that the right hand side is ${\mathcal
C}^1$ and hence the scalar curvature function is ${\mathcal C}^2$.

This then implies that any $k$-th order partial derivatives of the
functions $\psi^{o,j}$ and $\psi^{i,j}$, with $k \leq 4$, coincide
on $\del B_{1}$.

Therefore, we conclude that the function $\psi$ defined by $\psi : =
\psi^{o,j}$ in $\bar B_{2} -B_{1}$ and $\psi : = \psi^{i,j}$ in
$\bar B_{1} - B_{1/2}$ is ${\mathcal C}^{4}$ in $\bar B_{2} -
B_{1/2}$ and is a solution of the nonlinear elliptic partial
differential equation
\[
{\bf s} \,\left( i \, \del\, \bar \del ( \tfrac{1}{2}  \, |v|^{2}
+ \psi ) \right) = f .
\]
where $f$ is defined by \[ -df  =  i \del \, \bar \del (
\tfrac{1}{2} \, |v|^2 +  \psi) (X, -)
\]
and hence is a nonlocal first order differential operator in $\psi$.
It then follows from elliptic regularity theory together with a
bootstrap argument that the function $\psi$ is in fact smooth.

Hence, by gluing the \K\ metrics $\omega_{\e, {\bf h},{\bf k}}$
defined on $\bar M_{r_\e}$ with the metrics $\e^{2} \, \eta_{\e,\hat
a_j , \nu_j, X_j, \tilde h_j, \tilde k_j}$ defined on $\bar
N_{R_\e/\hat a_j}$, we will produced a \K\ metric on $M_{\e}$ which
has constant scalar curvature. This will end the proof of
Theorem~\ref{mainthm}.

It remains to explain how to find the boundary data
\[
{\bf h} =(h_1, \ldots, h_n), \quad {\bf k} = (k_1, \ldots, k_n),
\quad {\bf \tilde h} = (\tilde h_1, \ldots, \tilde h_n) \qquad
\mbox{and} \qquad {\bf \tilde k} = (\tilde k_1, \ldots, \tilde k_n)
\]
as well as the parameters ${\bf \tilde a} =  (\tilde a_1, \ldots,
\tilde a_n)$.

We change the boundary data functions $h_j$ and $k_j$ into $h'_j$
and $k'_j$ defined by
\[
\begin{array}{rllll}
h'_j & : = & (\tilde a_j - a_j) \, r_\e^{4-2m} \, \e^{2m-2}
+ h_j \\[3mm]
k'_j & : = &  4 \, (m-2) ( a_j - \tilde a_j) \, \e^{2m-2} \,
r_\e^{4-2m} + k_j
\end{array}
\]
when $m\geq 3$ and
\[
\begin{array}{rllll}
h'_j & : = & h_j \\[3mm]
k'_j & : = &  4 \,  ( a_j - \tilde a_j) \, \e^{2} \, + k_j
\end{array}
\]
when $m=2$.  Recall that the functions $k_j$ are assumed to have
mean $0$ while the functions $k'_j$ are not assumed to satisfy such
a constraint anymore. The role of the scalar $\tilde a_j - a_j$ is
precisely to recover this lost degree of freedom in the assignment
of the boundary data.

If  $k$ is a constant function of $\del B_1$, we extend the
definition of $W^{o}_{h,k}$ by setting
\[
W^o_{0 , k} : = \tfrac{k}{4 (m-2)} \, (|z|^{2-2m} - |z|^{4-2m})
\]
when $m\geq 3$ and by
\[
W^o_{0 , k} = \tfrac{k}{4} \, \log |z|^{2}
\]
when $m=2$.

We also do not assume anymore that $\tilde h_j$ and $\tilde k_j$
satisfy (\ref{eq:concon}) anymore. If  $h$ is a constant function of
$\del B_1$, we extend the definition of $W^{i}_{h,k}$ by setting
\[
W^i_{h , 0} : =  h
\]
Finally, we set
\[
\tilde h_j' := \e^{2} \, \tilde h_j \qquad \qquad \tilde k_j' :=
\e^{2} \, \tilde k_j
\]

With these new variables, the expansions for both $\psi^{o,j}$ and
$\psi^{i,j}$ can now be written as
\[
\begin{array}{rllll}
\psi^{o,j} & = & -  \tilde a_j \, r_\e^{4-2m} \,  \e^{2m-2} \,
|\cdot |^{4-2m} + W^o_{h_j', k_j'} + {\mathcal O}_{{\mathcal C}^{4, \alpha}} (r_\e^{4}) \\[3mm]
\psi^{i,j} & = & - \tilde a_j \, r_\e^{4-2m} \, \e^{2m-2} \, |\cdot
|^{4-2m} + \, W^i_{\tilde h'_j , \tilde k'_j }  + {\mathcal
O}_{{\mathcal C}^{4, \alpha}} (r_\e^{4}) .
\end{array}
\]
when $m\geq 3$ and
\[
\begin{array}{rllll}
\psi^{o,j} & = & \tilde a_j \,  \e^2 \, \log |\cdot | + W^o_{h_j', k_j'} + {\mathcal O}_{{\mathcal C}^{4, \alpha}} (r_\e^{4}) \\[3mm]
\psi^{i,j} & = & \tilde a_j \,  \e^2 \, \log |\cdot | + \,
W^i_{\tilde h'_j , \tilde k'_j } + {\mathcal O}_{{\mathcal C}^{4,
\alpha}} (r_\e^{4}) .
\end{array}
\]
when $m=2$. Observe that, since $\tilde h_j$ and $\tilde k_j$ a re
not assumed to satisfy (\ref{eq:concon}) anymore, this has changed
the value of $\psi^{i,j}$ by some constant which will not be
relevant for the computation of the corresponding \K\ form.

The system (\ref{eq:6.300}) we have to solve can now be written as~:
For all $j=1, \ldots, n$
\begin{equation}
\begin{array}{cccccccrllllll}
W^{o}_{h_j', k_j'} & = & W^{i}_{\tilde h_j', \tilde k_j'}  & +&
{\mathcal O}_{{\mathcal C}^{4, \alpha}} (r_\e^{4}) \\[3mm]
\del_r W^{o}_{h_j', k_j'} & = & \del_r W^{i}_{\tilde h_j', \tilde
k_j'} & + & {\mathcal O}_{{\mathcal C}^{3, \alpha}} (r_\e^{4}) \\[3mm]
\Delta W^{o}_{h_j', k_j'} & = & \Delta W^{i}_{\tilde h_j', \tilde
k_j'} & + & {\mathcal O}_{{\mathcal C}^{2, \alpha}}
(r_\e^{4}) \\[3mm]
\del_r \Delta W^{o}_{h_j', k_j'} & =  & \del_r  \Delta W^{i}_{\tilde
h_j', \tilde k_j'} & + & {\mathcal O}_{{\mathcal C}^{1, \alpha}}
(r_\e^{4}) \end{array} \label{eq:6.3000000}
\end{equation}
on $\del B_{1}$.

By definition of $W^o_{h,k}$ and $W^i_{h,k}$, the first equations
and third equations reduce to
\begin{equation}
\begin{array}{rllllll}
h_j' & = & \tilde h_j' + {\mathcal
O}_{{\mathcal C}^{4, \alpha}} (r_\e^{4}) \\[3mm]
k_j'  & = & \tilde k_j' + {\mathcal O}_{{\mathcal C}^{2, \alpha}}
(r_\e^{4}) \end{array} \label{eq:0077}
\end{equation}
Inserting these into the second and third sets of equations and
using the linearity of the mapping $(h,k) \longmapsto W^o_{h,k}$ and
$(h,k) \longmapsto W^i_{h,k}$, the second and third equations become
\begin{equation}
\begin{array}{cccccccrllllll}
\del_r W^{o}_{h_j', k_j'} & = & \del_r W^{i}_{h_j', k_j'} & + &
{\mathcal O}_{{\mathcal C}^{3, \alpha}} (r_\e^{4}) \\[3mm]
\del_r \Delta W^{o}_{h_j', k_j'} & =  & \del_r  \Delta W^{i}_{h_j',
k_j'} & + & {\mathcal O}_{{\mathcal C}^{1, \alpha}} (r_\e^{4})
\end{array} \label{eq:007}
\end{equation}
for all $j=1, \ldots, n$. We now make use of the following result
whose proof can be found in \cite{ap}~:
\begin{lemma}
The mapping
\[
\begin{array}{rclclll}
\mathcal P :& {\mathcal C}^{4,\alpha}(\del B_1 ) \times {\mathcal
C}^{2,\alpha} (\del B_1) & \longrightarrow & {\mathcal
C}^{3,\alpha}(\del B_1) \times {\mathcal C}^{1,\alpha}(\del B_1) \\[3mm]
& (h,k)  &\longmapsto    & (\partial_{r} \, (W^{i}_{h, k}- W^{o}_{h,
k}), \partial_{r} \, \Delta \, (W^{i}_{h, k}- W^{o}_{h, k})) ,
\end{array}
\]
is an isomorphism. \label{le:6.3}
\end{lemma}

Using Lemma~\ref{le:6.3}, (\ref{eq:007}) reduces to
\begin{equation}
\begin{array}{cccccccrllllll}
h_j'& = &  {\mathcal O}_{{\mathcal C}^{4, \alpha}} (r_\e^{4}) \\[3mm]
k_j'& =  & {\mathcal O}_{{\mathcal C}^{2, \alpha}} (r_\e^{4})
\end{array}
\label{008}
\end{equation}
for all $j=1, \ldots, n$. This, together with (\ref{eq:0077}),
yields a fixed point problem which can be written as
\[
({\bf h'} ,  {\bf \tilde h'} , {\bf  k'}, {\bf \tilde k'} ) = S_\e (
{\bf h'} , {\bf \tilde h'} , {\bf k'} , {\bf \tilde k} ) ,
\]
and we know from (\ref{eq:0077}) and (\ref{008}) that the nonlinear
operator $S_\e$ satisfies
\[
\| S_\e ({\bf h'} ,  {\bf \tilde h'} , {\bf  k'}, {\bf \tilde k'}  )
\|_{({\mathcal C}^{4, \alpha})^{2n} \times ({\mathcal C}^{2, \alpha}
)^{2n}} \leq c_0 \, r_\e^{4} ,
\]
for some constant $c_0 >0$ which does not depend on $\kappa$,
provided $\e \in (0, \e_\kappa)$. We finally choose
\[
\kappa = 2 \, c_0 ,
\]
and $\e \in (0, \e_{\kappa})$. We have therefore proved that $S_\e$
is a map from
\[
A_\e : = \left\{ ({\bf h'} ,  {\bf \tilde h'} , {\bf  k'}, {\bf
\tilde k'}  ) \in ({\mathcal C}^{4, \alpha})^{2n} \times ({\mathcal
C}^{2, \alpha})^{2n} \quad : \quad \| ({\bf h'} ,  {\bf \tilde h'} ,
{\bf  k'}, {\bf \tilde k'} ) \|_{({\mathcal C}^{4, \alpha})^{2n}
\times ({\mathcal C}^{2, \alpha})^{2n} } \leq \kappa \, r_\e^{4}
\right\} ,
\]
into itself. It follows from (\ref{eq:6.14}) and (\ref{eq:6.29})
that, reducing $\e_\kappa$ if this is necessary, $S_\e$ is a
contraction mapping from $A_\e$ into itself for all $\e \in (0,
\e_\kappa)$. Therefore, $S_\e$ has a fixed point in this set. This
completes the proof of the existence of a solution of
(\ref{eq:6.300}).

The proof of the existence on $M_{r_\e} $ of a \K\ form $\omega_\e$
which has constant scalar curvature is therefore complete. Observe
that the scalar curvature of $\omega$ and $\omega_\e$ are close
since the estimate
\[
|{\bf s} (\omega_\e) - {\bf s} (\omega)|\leq c \, \e^{2m-2}
\]
follows directly from the construction. We also have the estimates
\[
|\tilde a_j - a_j| \leq c \, r_\e^{2m}\, \e^{2m-2} = c \,
\e^{\tfrac{2}{2m+1}}
\]
which is the last estimate which appears in the statement of
Theorem~\ref{mainthm}

Since the construction of the \K\ form $\omega_\e$ is performed
using fixed point theorems for contraction mappings, it should be
clear that $\omega_\e$ depends continuously on the parameters of the
construction (such as the \K\ form $\omega$, the points $p_j$ and
the coefficients $a_j$. In particular, when $\frh'' =\{0\}$,
conditions (i), (ii) and (iii) in the statement of
Theorem~\ref{mainthm} is void and in constructing $\omega_\e$ one is
free to prescribe the parameters $a_j$. A simple degree argument
then shows that given $A \subset ({\mathbb R}^*_+)^n$ the image of
mapping
\[
(a_1, \ldots, a_n) \longmapsto (\tilde a_1, \ldots, \tilde a_n)
\]
contains $A$ provided $\e$ is chosen small enough. This completes
the proof of the last remark at the end of the statement of
Theorem~\ref{mainthm}. In the case where $\frh''\neq \{0\}$, one
needs to apply a modified version of the analysis of \cite{ap2}
which guaranties that the set of weights is open (keeping the
required symmetries) provided no nontrivial element of ${\frh''}$
vanishes at all points we blow up.

The proof of Proposition~\ref{mainprop} also follows from the
construction itself. Indeed, when $\omega$ is a constant scalar
curvature \K\ form, then $X_{\bf s} =0$. However, in the expansion
of Proposition~\ref{pr:6.22} one directly sees that the scalar
curvature of $\omega_\e$ will not be constant if the vector field
$Y'$ is not zero. Now $Y' =0$ if and only if $\sum_j a_j \,
\xi_{\omega} (p_j) =0$. This completes the proof of
Proposition~\ref{mainprop}

\section{Examples and Comments}

\subsection{Toric varieties}

If $M^m$ is a toric variety, then we can take $K=T^m$ and $\frh =
\frk$ in virtue of Proposition \ref{fixx}. 
Theorem \ref{mainthm-2} asserts that one can blow up {\em any
set of points contained in the fixed-point set of the torus-action}
and the weights $a_j >0$ can be chosen arbitrarily. This is because
the algebra of {\em good vector fields} that extends to the blown up
manifold is precisely the Lie algebra of the torus in this case, so
$\frh''=\{0\}$ and conditions (i) and (ii) become vacuous in this
case.

This type of examples leads naturally to some observations~:

\begin{enumerate}
\item[(i)]
As mentioned in the introduction, our procedure can be iterated
since blowing up a toric variety at such points preserves the toric
structure. Therefore, we obtain extremal metrics on any such
iterated blow up. Among this type of manifolds Donaldson \cite{do}
has studied one particular iterated blow up of ${\mathbb P}^{2}$, 
where all successive
blow ups take place at fixed points of the torus action at the previous
step. The number of iteration cannot be explicitly determined, yet
these manifolds fall into the category to which our result applies.
Nonetheless Donaldson analysis shows if on these manifolds
we take  \K\ classes sufficiently far from the boundary of the \K\ cone
no extremal representative exist.
This
shows that even for these deceptively simple examples the
understanding of the maximal range of application of our result
(namely the determination
of the optimal value of $\e_0$) is far from been trivial. \\[3mm]
\item[(ii)]
It is known that not every manifold admits an extremal metric but
Tian \cite{ti} has conjectured that every \K\ manifold degenerates
in a suitable sense to a manifold with an extremal metric.

\noindent We know two types of manifolds which do not carry any
extremal metric~: The projectivization of unstable rank two vector
bundles over compact Riemann surfaces \cite{bdb} (which verify
Tian's conjecture by their very construction) and some special
iterated blow up of ${\mathbb P}^{2}$ constructed by Levine
\cite{le}.

\noindent Our result can be used to shows how to fit these examples
in Tian's conjecture. For example, one observes that in Levine's
examples all the blow ups are of the type allowed by our
construction except the last one where two points which do not
correspond to vertices of the polytope are blown up. Nevertheless,
this last manifold degenerates on the manifold obtained by blowing
up one further vertex and then another point on the last exceptional
divisor and, by our main theorem extremal metrics do exist on this
limit manifold.
\end{enumerate}

Besides these general results, we now look at the effect of our construction
in specific cases.

Firstly, since ${\mathbb P}^{m}$ is a toric manifold the construction of the new
extremal metrics on its blow ups at $n\leq m+1$ linearly independent points
follows. The fact that the resulting extremal metrics have constant scalar curvature
iff we blow up $m+1$ points with equal weights follows from Proposition \ref{mainprop}
and from our preceding work \cite{ap2}. Nonetheless this last result can be easily
obtained directly.
Let us in fact look at the Futaki invariants of the resulting manifolds.
Let us recall Mabuchi's result \cite{ma1} stating that the Futaki invariant of a \K\ 
manifold vanishes if and only if the barycenter of the associated polytope lies 
in the origin.
It is in general a hard combinatorial task, knowing the polytope, to determine 
its barycenter
(as the reader of \cite{nak} can easily glance)
and this prevents us to state general results, yet restricting ourselves to 
projective spaces
we can get a clear picture.

The polytope associated to ${\mathbb P}^{m}$ with its \K -Einstein metric is
well known to be the simplex in $\R^m$ with vertices

\[
p_1=(1,0,\dots,0) \dots p_m=(0,\dots,0,1) \qquad p_{m+1}=(-1,\dots, -1)
\]

and the vertices are exactly the images of the fixed points of the torus action.
Let us recall also the effect of blowing up one of these points. For all the toric 
geometry we will use  
we refer to \cite{gu}.
Given a vertex $p_j$, the polytope associated to the manifold 
$(Bl_{p_j}M, \pi^*[\omega] - aPD[E])$ is the one
of M where the vertex $p_j$ is substituted by $m$ vertices $q_j^k$, $k=1,\dots m+1$, 
$k\neq j$, where $q_j^k = p_j + a(p_j - p_k)$.
Blowing up a vertex with weight $a>0$ has then the effect of {\em cutting} the starting
 polytope removing a simplex of  ``size" $a$.

In general it is possible to get that the barycenter stays in the origin even
blowing up fewer points than the whole set of vertices and with uncontrollable weights.
For example, let
us take as the base manifold the blow up of ${\mathbb P}^2$ at three
points which are not aligned. The polytope associated to the
canonical class $3\pi^*[\omega_{FS}] -  \left( PD[E_1] + PD[E_2] + PD
[E_3] \right)$ is the hexagon with vertices $(1,0)$, $(0,1)$,
$(1,1)$, $(-1,0)$, $(0,-1)$, and $(-1,-1)$. The existence of an
Einstein (constant scalar curvature) metric in this class is a
result of Siu \cite{si} and Tian-Yau \cite{ty}.  Blowing up two pairs
of vertices symmetric with respect to the origin (hence $4$ points)
with equal pairwise weights, we get new constant scalar curvature
metrics on the blow up.

On the other hand, in the case of ${\mathbb P}^{m}$, plugging in the above 
numbers, it is an elementary
calculation to get that the barycenter of the blow up is still in the origin iff  $n=m+1$ and 
$a_1 = \dots = a_{m+1}$. Calabi's result \cite{ca2}, stating that an extremal metric in 
a \K\ class whose Futaki invariant vanishes is of constant scalar curvature, 
gives a different proof of our characterization of \K\ constant scalar 
curvature metrics among  our new extremal ones.

Once we know that all manifolds obtained by blowing up any number $n$ of points
in general position in ${\mathbb P}^{m}$ admit an extremal metric 
(for $n\geq m+2$ these are in fact constant scalar curvature metrics \cite{ap2}), 
let us focus on the \K\ classes we can reach for $m=2$.

Recall that ${\mathbb P}^{1} \times {\mathbb P}^{1}$ is itself a toric variety whose polytope corresponding to the \K\ class $\alpha [\omega_{FS}] + \beta [\omega_{FS}]$,
$\alpha, \beta >0$, where $[\omega_{FS}^j$ is half the \K -Einstein metric on the $j$-th factor, is the rectangle with vertices $(\alpha, \beta), (-\alpha, \beta),(\alpha, -\beta), (-\alpha, -\beta)$. Let us also briefly recall the following classical construction~: take $p_1, p_2 \in {\mathbb P}^{2}$ and consider $M = Bl_{p_1,p_2}{\mathbb P}^{2}$.
$M$ contains three $(-1)$-curves, the two exceptional divisors $E_1$,$E_2$ and the proper transform $L$ of the line in ${\mathbb P}^{2}$ passing though $p_1$ and $p_2$.
$M$ is in fact biholomorphic to 
$Bl_q({\mathbb P}^{1} \times {\mathbb P}^{1})$ for some (hence {\em any}) choice of $q\in {\mathbb P}^{1} \times {\mathbb P}^{1}$. In fact, contracting (``blowing down") $L$
we get a manifold biholomorphic to ${\mathbb P}^{1} \times {\mathbb P}^{1}$ where the rulings correspond to the pencils of lines through $p_1$ and $p_2$.

Now, having called $A_1 = [{\mathbb P}^{1} \times \{pt\}]$, 
$A_2 = [\{pt\} \times {\mathbb P}^{1}]$, and $E$ the exceptional divisor in 
$Bl_q({\mathbb P}^{1} \times {\mathbb P}^{1})$, it is easy to check the 
correspondence of classes
\[
\alpha PD[A_1] + \beta PD[A_2] -\lambda PD[E] \leftrightarrow
(\alpha+\beta-\lambda)\pi^*[\omega_{FS}] - (\alpha-\lambda) PD[E_1]
-(\beta -\lambda)PD[E_2]\, .
\]

Hence our new extremal metrics on $Bl_q({\mathbb P}^{1} \times {\mathbb P}^{1})$,
which lie in the classes $\alpha PD[A_1]  + \beta PD[E_2]  - \e^2 \lambda PD[E]$, 
give extremal metrics in the classes of $Bl_{p_1,p_2}{\mathbb P}^{2}$
\[
\pi^*[\omega_{FS}] -  \, \tfrac{\alpha - \e^2 \lambda}{\alpha+\beta -\e^2}\, PD[E_1] -
\tfrac{\beta - \e^2\lambda}{\alpha+\beta-\e^2\lambda}\,  \, PD [E_2].
\]

Hence we have a whole neighborhood of the boundary line $a+b=1$ in the \K\ cone of
$Bl_{p_1,p_2}{\mathbb P}^{2}$,  where $\pi^*[\omega_{FS}] -  \, a PD[E_1] -
 b\,  \, PD [E_2]$, $a+b \leq 1$, $a,b >0$.
 This proves Corollary \ref{claproj}, $(1)$.

 Another similar construction, known as Cremona transformation (see e.g. \cite{har} pages 397-399), allows us to prove Corollary \ref{claproj}, $(2)$. We wish to thank M.Abreu for bringing it to our attention.
 
 This time we construct an automorphism of $Bl_{p_1,p_2,p_3}{\mathbb P}^{2}$, 
when the points do not lie on a  complex line, in the following way~:
call $l_{jk}$ the lines in ${\mathbb P}^{2}$ through $p_j$ and $p_k$ and $L_{jk}$ their proper transforms. The key observation this time is that blowing down $L_{jk}$ 
we are left with a new copy of ${\mathbb P}^{2}$, where the new coordinate lines 
are the exceptional divisors of the original blow up. The resulting automorphism of 
$Bl_{p_1,p_2,p_3}{\mathbb P}^{2}$ has the following action in cohomology, where we indicate by $F_j$ the exceptional divisors in the new copy of 
$Bl_{p_1,p_2,p_3}{\mathbb P}^{2}$~:
\[
\pi^*[\omega_{FS}] - PD[E_1] - PD[E_2] \leftrightarrow PD[F_1] ,
\]
\[
\pi^*[\omega_{FS}] - PD[E_1] - PD[E_3] \leftrightarrow PD[F_2] ,
\]
\[
\pi^*[\omega_{FS}] - PD[E_2] - PD[E_3] \leftrightarrow PD[F_3] ,
\]
\[
\pi^*[\omega_{FS}]  \leftrightarrow 2\pi^*[\omega_{FS}]  - PD[F_1] - PD[F_2] - PD[F_3].
\]

Hence
\[
\pi^*[\omega_{FS}] -a  PD[E_1] - b PD[E_2] - c PD[E_3] 
\]
corresponds to
\[
(2-a-b-c)\pi^*[\omega_{FS}] - (1-a-b) PD[F_1] - (1-a-c) PD[F_2] - (1-b-c) PD[F_3]
\]

This gives the seeked \K\ classes as claimed.

\subsection{Sporadic examples}

In the case where the manifold we study is the projective space, we
can also study the existence existence of extremal metrics on the
blown up manifold when the position of the blown up points is not
generic, hence leaving the toric world. 
For example in \cite{ap2} it was shown how to construct \K\
constant scalar curvature metrics on ${\mathbb P}^2$ blown up at $4$
points, $3$ of which were lying on a line. Considering extremal
metrics instead on constant scalar curvature metrics, one can do
more~: for example, consider the points
\[
p_1= [0 : \tfrac{1}{\sqrt{2}} :  \tfrac{1}{\sqrt{2}}], \, \,
p_2 = [0, \tfrac{\alpha}{\sqrt{|\alpha|^2 + |\beta|^2}} : \tfrac{\beta}{\sqrt{|\alpha|^2 + |\beta|^2}}] \,\, 
p_3 = [0 : \tfrac{\beta}{\sqrt{|\alpha|^2 + |\beta|^2}} : \tfrac{\alpha}{\sqrt{|\alpha|^2 + |\beta|^2}}]
\]
and the group $K = S^{1}$ whose action on ${\mathbb P}^2$ is given
by
\[
\begin{array}{crcllllll}
S^1 \times {\mathbb P}^2 & \longrightarrow & {\mathbb P}^2 \\[3mm]
(\theta , [z^{1} : z^{2}: z^{3}]) & \longmapsto & [\theta^{-2} \,
z^{1} : \theta \, z^{2} : \theta \, z^{3} ]
\end{array}
\]

Of course $p_1, p_2$ and $p_3$ are fixed under the action of $K$,
but we want to impose more symmetries working equivariantly with
respect to a discrete group $A$ of permutations of the last two
coordinates. It is easy to check that the space of vector fields of
$\frh$ which are invariant under the action of $A$ is given by
\[
\frh _{A} = \mbox{Span} \{ \Re(z^{2} \, \partial _{z^{3}} + z^{3} \,
\partial _{z^{2}}) ,  \Re(2 \, z^{1} \, \partial_{z^{1}}
 - z^{2} \, \partial_{z^{2}} - z^{3} \, \partial_{z^{3}}) \}
\]
It is immediate to see that
\[ \frh' _{A} = \mbox{Span} \{\Re(2 \, z^{1} \, \partial_{z^{1}}
 - z^{2} \, \partial_{z^{2}} - z^{3} \, \partial_{z^{3}} )\}
 \]
Observe that all points belong to the zero set of the vector fields
in $\frh'$. Condition (ii) in Theorem~\ref{mainthm} is fulfilled
provided
\[
\Re \, (\alpha \, \bar \beta) < 0
\]
with $a_1 = - \tfrac{2 \Re \, (\alpha \, \bar \beta)}{|\alpha|^2 + |\beta|^2}$ and 
$a_2=a_3 =1$, while condition
(iii) holds since $ \Re(z^{2} \, \partial _{z^{3}} + z^{3} \,
\partial _{z^{2}}) $ does not vanish at $p_2$ and at $p_3$. 
It is important to observe that $a_1$ automatically satisfies the upper bound $a_1\leq 1$.
In fact, it has been shown by A. Della Vedova \cite{dv} that for $a_1 >2$
the corresponding polarized manifold is not relative K-stable, hence forbidding the existence of an extremal metric in the corresponding \K\ class
thanks to Szekelehidi's result \cite{sz}.

This
gives the construction of extremal \K\ metrics on the blow up
${\mathbb P}^2$ at $p_1, p_2$ and $p_3$.
The fact that the corresponding
metrics do not have constant scalar curvature follows directly from Proposition~\ref{mainprop},
or implicitely observing that the resulting manifold does not satisfy the Matsushima-Licherovitz
obstruction so in fact it does not admit constant scalar curvature metrics in {\em any} \K\ class.

The same remark we observed after Corollary \ref{coproj} also holds
for this example. Indeed, the position of three (even aligned)
points in ${\mathbb P}^2$ can be {\em a priori} prescribed (just
leaving them on some line) with the use of an appropriate
automorphism. Therefore the above calculation implies that {\em any}
set of three aligned points can be blow up and extremal metrics on
the blown up manifold can be found even if in the initial
coordinates the symmetries of the above example are not present. Of
course the change of coordinates required to put the initial set of
points into the above one might change the base Fubini-Study metric
we work with, but this does not change its \K\ class. This
discussion can be summarized in the~:
\begin{corol}
Given three aligned points $p_1, p_2, p_3$ in ${\mathbb P}^2$ and
weights $a_1, a_2, a_3$ satisfying $a_j=a_k$ for some $j,k$ and 
$a_l \leq a_j$ for all $l,j$, there
exists $\e_0 >0$ and for all $\e \in (0, \e_0)$ there exists an
extremal \K\ form of non constant scalar curvature $\omega_\e$ on
the blow up of ${\mathbb P}^2$ at $p_1, p_2, p_3$ with
\[
\omega_\e \in  \pi^*[\omega_{FS}] -  \e^2 \, \left( a_1\, PD[E_1] +
a_2\, PD[E_2] +  a_3 \, PD [E_3] \right)
\]
\label{co:align}
\end{corol}

\noindent As expected, adding points to be blown up, also on the
same line, makes things even simpler. For example let us work out
the situation where $4$ aligned points are to be blown up. In this
case we can avoid using extra symmetries and we can work directly
with a connected group of isometries. Therefore, we now consider the
points
\[
p_1 = [0: 1 : 0], \quad p_2=[0: \tfrac{1}{\sqrt{3}} : \tfrac{1+i}{\sqrt{3}}], \qquad p_3 = [0 : \tfrac{1}{\sqrt{5}} : \tfrac{-2-i}{\sqrt{5}}],
\qquad \mbox{and} \qquad p_4=[0 : \tfrac{1}{\sqrt{2}} : \tfrac{1}{\sqrt{2}}]
\]
and the group $K = S^{1}$, whose action on ${\mathbb P}^2$ given by
\[
\begin{array}{crcllllll}
S^1 \times {\mathbb P}^2 & \longrightarrow & {\mathbb P}^2 \\[3mm]
(\alpha , [z^{1},z^{2},z^{3}]) & \longmapsto & [\alpha  \, z^{1},
z^{2}, z^{3} ]
\end{array}
\]

Of course $p_1, \dots, p_4$ are fixed by the action of $K$. It is
easy to check that $\frh$ is now given by
\[
\frh = \mbox{Span} \{ \Re(z^{2} \partial _{z^{3}} + z^{3} \partial
_{z^{2}}) ,\Re( i \, (z^{2} \partial _{z^{3}} - z^{3} \partial _{z^{2}}))
, \Re(2 z^{1} \partial_{z^{1}} - z^{2}\partial_{z^{2}} -
z^{3}\partial_{z^{3}}), \Re( z^{2}\partial_{z^{2}} - z^{3}\partial_{z^{3}})
\}
\]
We choose
\[
\frh' = \mbox{Span} \{ \Re(2 z^{1} \partial_{z^{1}}
 - z^{2}\partial_{z^{2}} - z^{3}\partial_{z^{3}} )\}
 \]
Observe that all points belong to the zero set of the vector fields
in $\frh'$ and that none of the nontrivial elements
of $\frh''$ vanish at all the $p_j$. Our construction then gives extremal (non constant
scalar curvature) \K\ metrics on the blow up of ${\mathbb P}^2$ at
$p_1, \dots, p_4$ with weights $a_1 =1$, $a_2 = 3$, $a_3 =  5$
and $a_4= 2$.

\noindent These examples can be easily extended to projective spaces
of any dimension.

\end{document}